# Artificial viscosity to cure the shock instability in high-order Godunov-type schemes

Alexander V. Rodionov

*Russian Federal Nuclear Center – All-Russian Scientific Research Institute of Experimental Physics (RFNC-VNIIEF), Sarov, Nizhny Novgorod region, 607190, Russia*

*E-mail address*: avrodionov@rambler.ru

ABSTRACT. The artificial viscosity approach for curing the carbuncle phenomenon (a numerical problem, also known as the shock instability) in shock-capturing methods has been recently presented and successfully tested on the first-order schemes in two- and three-dimensional simulations. The present study extends the proposed approach to the case of using high-order Godunov-type schemes. Several implementations of well-known schemes were selected for the study. They involve the MUSCL and WENO data reconstructions in space along with the Runge-Kutta and Hancock-type time stepping techniques. Numerous computations of the Quirk-type test problems and other popular tests were performed to examine and tune the artificial viscosity approach as applied to the selected schemes. As a result of this study (1) the principal coefficient in the artificial viscosity model was adjusted and (2) some methodological suggestions for substantial weakening of the post-shock oscillations were stated.

## 1. Introduction

At present, computational fluid dynamics makes extensive use of shock-capturing schemes (methods) based on the exact or approximate solution of the Riemann problem; such schemes are called Godunov-type or (in a broader sense) upwind schemes. Offering a number of advantages, they also have drawbacks, of which the carbuncle phenomenon is the most severe [1, 2]. Over the thirty years gone from the time the phenomenon was noticed and described, a large number of papers have been published seeking to study and cure this flaw. (Note that the term *carbuncle phenomenon* is often used interchangeably with the terms *carbuncle instability* or *shock instability*.)

In [3] the author presented a summary on the carbuncle phenomenon and proposed a new approach for curing the problem. Its idea is to introduce some dissipation in the form of the right-hand side of Navier-Stokes equations into the basic method of solving Euler equations; in so doing, the molecular viscosity coefficient is replaced by the artificial viscosity coefficient. The new cure for the carbuncle flaw was tested and tuned for the case of using first-order schemes in two-dimensional simulations. Its efficiency was demonstrated on several well-known test problems. In [4] the artificial viscosity approach was extended to the case of three-dimensional simulations.

This work tests the applicability of the new cure for the shock instability in high-order Godunov-type schemes. The paper is organized as follows. The next section (section 2) describes the schemes selected for the study. Section 3 documents their detailed testing, as a result of which some methodological suggestions are stated and the principal coefficient in the artificial viscosity model is adjusted. Section 4 demonstrates the efficiency of the artificial viscosity approach on various test problems.

## 2. Numerical schemes

A variety of schemes, which can be classified as high-order Godunov-type schemes, have been developed to date. As distinct from first-order schemes, they differ not only in their Riemann solvers, but also in a variety of techniques for improving their accuracy in space and time. Moreover, the effect of applying a specific high-order scheme to some test problem can also depend on some minor (nonprincipal) features of its realization. For this reason, it does not seem possible to test the efficiency of the proposed approach for curing the carbuncle phenomenon on the whole variety of this class of schemes.

Instead, we restrict this work to considering implementations of some well-known schemes. As a result of their comprehensive testing we are going to develop general recommendations, which will definitely prove to be useful also for other schemes in the class of interest (high-order Godunov-type schemes). The schemes selected for testing are described below in this section.

### 2.1. The Euler equations

We describe the schemes used in this work in the context of solving two-dimensional equations of gas dynamics, which in Cartesian coordinates $xy$ read

$$\frac{\partial \mathbf{U}}{\partial t} + \frac{\partial \mathbf{F}_x}{\partial x} + \frac{\partial \mathbf{F}_y}{\partial y} = \text{RHS}, \qquad \mathbf{U} = \begin{bmatrix} \rho \\ \rho u_x \\ \rho u_y \\ \rho h_0 - p \end{bmatrix}, \quad \mathbf{F}_x = \begin{bmatrix} \rho u_x \\ \rho u_x^2 + p \\ \rho u_x u_y \\ \rho u_x h_0 \end{bmatrix}, \quad \mathbf{F}_y = \begin{bmatrix} \rho u_y \\ \rho u_x u_y \\ \rho u_y^2 + p \\ \rho u_y h_0 \end{bmatrix}, \tag{1}$$

where RHS = 0 for the Euler equations, $\mathbf{u} = (u_x, u_y)$ are the velocity vector and its components, $\rho$ is the density, $p$ is the pressure and $h_0$ is the specific total enthalpy, which for a polytropic gas is given by

$$h_0 = \frac{1}{2}\left(u_x^2 + u_y^2\right) + \frac{\gamma p}{(\gamma - 1)\rho},$$

with $\gamma$ denoting the ratio of specific heats.

### 2.2 Space-time discretization

We apply the finite volume approach to solve (1) numerically. In so doing we consider a sufficiently smooth structured grid and introduce curvilinear coordinates $\xi\eta$ that transform the grid in physical space $xy$ to a rectangular grid in computational space with the grid spacing $\Delta\xi = \Delta\eta = 1$. The following geometric parameters of the grid will take the place of metrics (hereinafter the grid indices $i$ and $j$ correspond to the coordinates $\xi$ and $\eta$): $V_{i,j}$, the cell volume; $(\mathbf{S}_\xi)_{i+\frac{1}{2},j}$ and $(\mathbf{S}_\eta)_{i,j+\frac{1}{2}}$, the area vectors of the cell faces in two grid directions (each area vector points in the direction of increase in the associated coordinate).

Assume that the solution at a given time level $n$ is defined, meaning that for the time $t = t^n$ we know the average values $\mathbf{U}_{i,j}^n$ within each cell. In order to update the solution to the next time level $n+1$ we integrate (1) over the cell and the time between $t^n$ and $t^{n+1} \equiv t^n + \Delta t$. The result is

$$\mathbf{U}_{i,j}^{n+1} = \mathbf{U}_{i,j}^n - \frac{\Delta t}{V_{i,j}}\Big[\left(\mathbf{F}_\xi\right)_{i+\frac{1}{2},j} - \left(\mathbf{F}_\xi\right)_{i-\frac{1}{2},j} + \left(\mathbf{F}_\eta\right)_{i,j+\frac{1}{2}} - \left(\mathbf{F}_\eta\right)_{i,j-\frac{1}{2}}\Big], \tag{2}$$

where $\mathbf{F}_\xi = S_{\xi x}\mathbf{F}_x + S_{\xi y}\mathbf{F}_y$ and $\mathbf{F}_\eta = S_{\eta x}\mathbf{F}_x + S_{\eta y}\mathbf{F}_y$ are the flux vectors across the cell faces in two grid directions.

The way of computing the fluxes depends on the specific scheme. In the first-order Godunov-type schemes we assume that the flow variables within each computing cell are constant, and the fluxes are computed from one of the Riemann solvers, the exact or approximate solution of the Riemann problem. Any of these solutions can be expressed as

$$\left(\mathbf{F}_\xi\right)_{i+\frac{1}{2},j} = \mathbf{F}^{RS}\left(\mathbf{Q}_{i,j}^n, \mathbf{Q}_{i+1,j}^n, \left(\mathbf{S}_\xi\right)_{i+\frac{1}{2},j}\right), \qquad \left(\mathbf{F}_\eta\right)_{i,j+\frac{1}{2}} = \mathbf{F}^{RS}\left(\mathbf{Q}_{i,j}^n, \mathbf{Q}_{i,j+1}^n, \left(\mathbf{S}_\eta\right)_{i,j+\frac{1}{2}}\right),$$

where $\mathbf{Q} \equiv \left(u_x, u_y, \rho, p\right)^T$ and $\mathbf{Q}_{i,j}^n = \mathbf{Q}\left(\mathbf{U}_{i,j}^n\right)$.

### 2.3. Reconstructions

Increasing the accuracy of the Godunov-type schemes is usually obtained by data reconstruction. In the second-order MUSCL-type schemes, the piecewise linear distribution of data is used instead of the piecewise constant distribution. The distributions of flow variables are reconstructed based on the known values of $\mathbf{Q}_{i,j}^n$ positioned at the cell centers. On a sufficiently smooth structured grid the reconstruction process reduces to calculating the so-called slopes, $\Delta_\xi \mathbf{Q}_{i,j}^n$ and $\Delta_\eta \mathbf{Q}_{i,j}^n$. The boundary extrapolated values are then found as

$$\mathbf{Q}_{i+\frac{1}{2}-,j}^n = \mathbf{Q}_{i,j}^n + \frac{1}{2}\Delta_\xi \mathbf{Q}_{i,j}^n, \qquad \mathbf{Q}_{i-\frac{1}{2}+,j}^n = \mathbf{Q}_{i,j}^n - \frac{1}{2}\Delta_\xi \mathbf{Q}_{i,j}^n,$$

$$\mathbf{Q}_{i,j+\frac{1}{2}-}^n = \mathbf{Q}_{i,j}^n + \frac{1}{2}\Delta_\eta \mathbf{Q}_{i,j}^n, \qquad \mathbf{Q}_{i,j-\frac{1}{2}+}^n = \mathbf{Q}_{i,j}^n - \frac{1}{2}\Delta_\eta \mathbf{Q}_{i,j}^n,$$

and for the fluxes we can write

$$\left(\mathbf{F}_\xi\right)_{i+\frac{1}{2},j} = \mathbf{F}^{RS}\left(\mathbf{Q}_{i+\frac{1}{2}-,j}^n, \mathbf{Q}_{i+\frac{1}{2}+,j}^n, \left(\mathbf{S}_\xi\right)_{i+\frac{1}{2},j}\right), \qquad \left(\mathbf{F}_\eta\right)_{i,j+\frac{1}{2}} = \mathbf{F}^{RS}\left(\mathbf{Q}_{i,j+\frac{1}{2}-}^n, \mathbf{Q}_{i,j+\frac{1}{2}+}^n, \left(\mathbf{S}_\eta\right)_{i,j+\frac{1}{2}}\right).$$

Note that reconstructing the piecewise linear distribution of $\mathbf{Q}$ (instead of reconstructing the conservative variables) is a standard practice for MUSCL-type schemes (see e.g. the review article by van Leer [5]).

The slopes are computed by the same slope-limiter algorithm for each grid direction independently (one-dimensional reconstruction). Below we consider three slope limiters: minmod, vanLeer and MC (the limiters are named in accordance with the book [6]; MC = monotonized central-difference limiter). They can be written as

$$\Delta f_i = \text{minmod}\left(\Delta f_{i-\frac{1}{2}}, \Delta f_{i+\frac{1}{2}}\right), \quad \text{minmod}\left(a,b\right) \equiv \begin{cases} a, & \text{if } |a| \le |b| \text{ and } ab > 0, \\ b, & \text{if } |b| < |a| \text{ and } ab > 0, \\ 0, & \text{if } ab \le 0, \end{cases}$$

$$\Delta f_i = \text{vanLeer}\left(\Delta f_{i-\frac{1}{2}}, \Delta f_{i+\frac{1}{2}}\right), \qquad \text{vanLeer}\left(a,b\right) \equiv \begin{cases} \dfrac{2ab}{a+b}, & \text{if } ab > 0, \\ 0, & \text{if } ab \le 0, \end{cases} \qquad (3)$$

$$\Delta f_i = \text{MC}\left(\Delta f_{i-\frac{1}{2}}, \Delta f_{i+\frac{1}{2}}\right), \qquad \text{MC}\left(a,b\right) \equiv \text{minmod}\left(\dfrac{a+b}{2}, \, 2\,\text{minmod}\left(a,b\right)\right).$$

where $\Delta f_{i-\frac{1}{2}} = f_i - f_{i-1}$ and $\Delta f_{i+\frac{1}{2}} = f_{i+1} - f_i$.

The reconstructions based on the limiters (3) belong to the class of three-point reconstructions satisfying the TVD condition. Of these, the minmod limiter is the most dissipative.

The effect of the reconstruction depends not only on the chosen slope limiter, but also on the variables, based on which the reconstruction is performed. The most straightforward way is to perform a reconstruction based on *primitive variables*, i.e. components of the vector $\mathbf{Q}$. This way, however, often leads to poor results (see an example in section 3.1). It is more consistent to perform a reconstruction based on *characteristic variables*. In the present work, this approach is implemented as follows.

Suppose we need to perform a reconstruction in cell $(i, j)$ in grid direction $\xi$. Two sets of primitive slopes, $\Delta_\xi \mathbf{Q}_{i-\frac{1}{2},j}^n = \mathbf{Q}_{i,j}^n - \mathbf{Q}_{i-1,j}^n$ and $\Delta_\xi \mathbf{Q}_{i+\frac{1}{2},j}^n = \mathbf{Q}_{i+1,j}^n - \mathbf{Q}_{i,j}^n$, are first calculated, and then they are converted into characteristic slopes, $\Delta_\xi \mathbf{Z}_{i-\frac{1}{2},j}^n$ and $\Delta_\xi \mathbf{Z}_{i+\frac{1}{2},j}^n$, using the transformation matrix $\mathbf{A} = \left(\partial \mathbf{Z} / \partial \mathbf{Q}\right)_{i,j}^n$. Next, one of the limiters (3) is applied to the characteristic slopes and the resulting vector $\Delta_\xi \mathbf{Z}_{i,j}^n$ is converted into $\Delta_\xi \mathbf{Q}_{i,j}^n$ by the inverse transformation matrix $\mathbf{A}^{-1}$; this completes the reconstruction process. The two transformation matrices are given by

$$\mathbf{A} = \begin{bmatrix} n_x & n_y & 0 & 1/\rho c \\ n_x & n_y & 0 & -1/\rho c \\ 0 & 0 & 1 & -1/c^2 \\ n_y & -n_x & 0 & 0 \end{bmatrix}, \qquad \mathbf{A}^{-1} = \begin{bmatrix} n_x/2 & n_x/2 & 0 & n_y \\ n_y/2 & n_y/2 & 0 & -n_x \\ \rho/2c & -\rho/2c & 1 & 0 \\ \rho c/2 & -\rho c/2 & 0 & 0 \end{bmatrix},$$

where $\rho = \rho_{i,j}^n$, $c = \sqrt{\gamma p_{i,j}^n / \rho_{i,j}^n}$, $\mathbf{n} \equiv (n_x, n_y) = \mathbf{S}_\xi / \|\mathbf{S}_\xi\|$, $\mathbf{S}_\xi = \left[ (\mathbf{S}_\xi)_{i-\frac{1}{2},j} + (\mathbf{S}_\xi)_{i+\frac{1}{2},j} \right] / 2$.

In addition to the piecewise linear reconstructions, in this work we employ the fifth-order reconstruction WENO5. A detailed description of this reconstruction is provided in [7], and its FORTRAN code is presented in [8]. We have tested alternative ways of using the characteristic variables in the WENO5 reconstruction: (1) applying the matrices of left and right eigenvectors of the Jacobian $\partial \mathbf{F}_\xi / \partial \mathbf{U}$ at the cell interface ($i+\frac{1}{2}$, $j$) (frozen values) when reconstructing $\mathbf{U}_{i+\frac{1}{2}-,j}$ and $\mathbf{U}_{i+\frac{1}{2}+,j}$ (as in [7]), and (2) applying these matrices at the cell center ($i$, $j$) when reconstructing $\mathbf{U}_{i-\frac{1}{2}+,j}$ and $\mathbf{U}_{i+\frac{1}{2}-,j}$ (as in [8]). The artificial viscosity approach was found to be equally effective in both cases. In what follows we present only the data obtained with the first version of using the characteristic variables in the WENO5 reconstruction.

### 2.4 Runge–Kutta time stepping

Rewrite Eq. (2) in a compact form:

$$\mathbf{U}_{i,j}^{n+1} = \mathbf{U}_{i,j}^n + \Delta t \mathbf{L}_{i,j}(\mathbf{Q}^n), \qquad (4)$$

where the spatial operator $\mathbf{L}_{i,j}$ expresses the total flux across all the cell faces. Its calculation implies a data reconstruction in each grid direction and solution of the Riemann problem on each cell face.

Equation (4) presents the forward Euler time discretization method which is first-order accurate. To increase the accuracy in time one commonly applies Runge–Kutta methods. In this work we use second- and third-order Runge-Kutta methods [9] (RK2 and RK3, respectively) having the following form

- RK2 method

$$\mathbf{U}_{i,j}^{(1)} = \mathbf{U}_{i,j}^n + \Delta t \mathbf{L}_{i,j}(\mathbf{Q}^n),$$

$$\mathbf{U}_{i,j}^{n+1} = \frac{1}{2}\mathbf{U}_{i,j}^n + \frac{1}{2}\mathbf{U}_{i,j}^{(1)} + \frac{1}{2}\Delta t \mathbf{L}_{i,j}(\mathbf{Q}^{(1)}).$$

- RK3 method

$$\mathbf{U}_{i,j}^{(1)} = \mathbf{U}_{i,j}^n + \Delta t \mathbf{L}_{i,j}(\mathbf{Q}^n),$$

$$\mathbf{U}_{i,j}^{(2)} = \frac{3}{4}\mathbf{U}_{i,j}^n + \frac{1}{4}\mathbf{U}_{i,j}^{(1)} + \frac{1}{4}\Delta t \mathbf{L}_{i,j}(\mathbf{Q}^{(1)}),$$

$$\mathbf{U}_{i,j}^{n+1} = \frac{1}{3}\mathbf{U}_{i,j}^n + \frac{2}{3}\mathbf{U}_{i,j}^{(2)} + \frac{2}{3}\Delta t \mathbf{L}_{i,j}(\mathbf{Q}^{(2)}).$$

Note that compared to the Euler method, the RK2 and RK3 methods are respectively two and three times more time-consuming.

### 2.5. Predictor-corrector technique

To achieve the second order of accuracy in time, other methods can be used, including the predictor-corrector technique described in [10, 11]. This specific technique is a kind of the Hancock scheme [12, 13], so in this work we call it the Hancock-Rodionov scheme (or HR scheme).

The predictor-corrector procedure of the HR scheme can be written in its compact form as

$$\tilde{\mathbf{U}}_{i,j}^{n+1} = \mathbf{U}_{i,j}^n + \Delta t \mathbf{L}_{i,j}^{pred}(\mathbf{Q}^n),$$

$$\mathbf{U}_{i,j}^{n+1} = \mathbf{U}_{i,j}^n + \Delta t \mathbf{L}_{i,j}^{corr}(\mathbf{Q}^{n+\frac{1}{2}}), \qquad \mathbf{Q}_{i,j}^{n+\frac{1}{2}} = \frac{1}{2}\left(\mathbf{Q}_{i,j}^n + \tilde{\mathbf{Q}}_{i,j}^{n+1}\right).$$

The spatial operators $\mathbf{L}_{i,j}^{pred}$ and $\mathbf{L}_{i,j}^{corr}$ are different here. As distinct from $\mathbf{L}_{i,j}$, the operator $\mathbf{L}_{i,j}^{pred}$ does not include solving the Riemann problem. Instead, the fluxes are calculated based on the boundary extrapolated values, which are interior with respect to the cell being integrated. Thus, at the predictor step, the fluxes for the cell $(i, j)$ are calculated by simple relations:

$$\left(\mathbf{F}_\xi\right)_{i-\frac{1}{2},j} = \mathbf{F}\left(\mathbf{Q}^n_{i-\frac{1}{2}+,j}, \left(\mathbf{S}_\xi\right)_{i-\frac{1}{2},j}\right), \qquad \left(\mathbf{F}_\xi\right)_{i+\frac{1}{2},j} = \mathbf{F}\left(\mathbf{Q}^n_{i+\frac{1}{2}-,j}, \left(\mathbf{S}_\xi\right)_{i+\frac{1}{2},j}\right),$$

$$\left(\mathbf{F}_\eta\right)_{i,j-\frac{1}{2}} = \mathbf{F}\left(\mathbf{Q}^n_{i,j-\frac{1}{2}+}, \left(\mathbf{S}_\eta\right)_{i,j-\frac{1}{2}}\right), \qquad \left(\mathbf{F}_\eta\right)_{i,j+\frac{1}{2}} = \mathbf{F}\left(\mathbf{Q}^n_{i,j+\frac{1}{2}-}, \left(\mathbf{S}_\eta\right)_{i,j+\frac{1}{2}}\right).$$

At the corrector step, we use the operator $\mathbf{L}^{corr}_{i,j}$, which includes solving the Riemann problem (like in the basic operator $\mathbf{L}_{i,j}$), but does not recompute the slopes: at an intermediate time level $n+\frac{1}{2}$, the same slopes are used as those obtained at the predictor step ($\varDelta_\xi\mathbf{Q}^n_{i,j}$ and $\varDelta_\eta\mathbf{Q}^n_{i,j}$).

Note that although the HR scheme requires less computing time than the RK2 method, it usually provides better accuracy.

### 2.6. The Navier-Stokes equations

The approach we employ to cure the carbuncle instability involves the Navier-Stokes equations. In this case the right-hand side in (1) reads

$$\text{RHS} = \frac{\partial \mathbf{F}^v_x}{\partial x} + \frac{\partial \mathbf{F}^v_y}{\partial y}, \quad \mathbf{F}^v_m = \left[0,\ \tau_{xm},\ \tau_{ym},\ q_m + u_x\tau_{xm} + u_y\tau_{ym}\right]^T, \quad m = x, y, \qquad (5)$$

where the stress tensor ($\tau_{nm}$) and the heat flux ($q_m$) can be written as

$$\tau_{nm} = \mu\left[\frac{\partial u_n}{\partial m} + \frac{\partial u_m}{\partial n} - \frac{2}{3}(\nabla\mathbf{u})\delta_{nm}\right], \qquad \nabla\mathbf{u} = \frac{\partial u_x}{\partial x} + \frac{\partial u_y}{\partial y},$$

$$q_m = \lambda\frac{\partial T}{\partial m} = \frac{\mu}{\text{Pr}}\left[\frac{\partial h_0}{\partial m} - \left(u_x\frac{\partial u_x}{\partial m} + u_y\frac{\partial u_y}{\partial m}\right)\right], \quad m, n = x, y,$$

with $\delta_{nm}$ denoting the Kronecker delta. In our approach the viscosity coefficient $\mu = \mu_{AV}$ and the Prandtl number $\text{Pr} = \frac{3}{4}$.

The space discretization of (5) takes the form

$$\text{RHS} = \frac{1}{V_{i,j}}\left[\left(\mathbf{F}^v_\xi\right)^n_{i+\frac{1}{2},j} - \left(\mathbf{F}^v_\xi\right)^n_{i-\frac{1}{2},j} + \left(\mathbf{F}^v_\eta\right)^n_{i,j+\frac{1}{2}} - \left(\mathbf{F}^v_\eta\right)^n_{i,j-\frac{1}{2}}\right],$$

where $\mathbf{F}^v_\xi = S_{\xi x}\mathbf{F}^v_x + S_{\xi y}\mathbf{F}^v_y$ and $\mathbf{F}^v_\eta = S_{\eta x}\mathbf{F}^v_x + S_{\eta y}\mathbf{F}^v_y$.

To compute viscous fluxes we use explicit central difference approximations, as described in [4]. The resulting value of the right-hand side of Navier-Stokes equations is introduced into the operator $\mathbf{L}_{i,j}$, or any of its variants, as source terms, which remain frozen (not recomputed) at all stages of marching the solution to the next time level $n+1$.

### 2.7. Artificial viscosity model

In high-order schemes we use the same artificial viscosity model as the one chosen for the first-order schemes [3]. This model defines the artificial viscosity coefficient as

$$\mu_{AV} = \begin{cases} C_{AV}\rho h^2\sqrt{(\nabla\mathbf{u})^2 - (C_{th}a/h)^2}, & \text{if } -(\nabla\mathbf{u}) > C_{th}a/h, \\ 0, & \text{otherwise}, \end{cases}$$

where $a$ is the sound speed, $h$ is the characteristic mesh size, $C_{AV}$ is the dimensionless parameter, $C_{th} = 0.05$ is the coefficient in the compression intensity threshold that restricts the effect of artificial viscosity to the shock layer only. In the two-dimensional case with non-square cells the characteristic mesh size is computed as $h = \max(d_1, d_2)/\sqrt{2}$, with $d_1$ and $d_2$ denoting the lengths of the cell diagonals. For rectangular cells in Cartesian coordinates it reduces to $h = \sqrt{\left(h_x^2 + h_y^2\right)/2}$ . In the three-dimensional case we use the expression for $h$ that was derived in [4].

Adaptation of this artificial viscosity model to high-order Godunov-type schemes will involve tuning the value of $C_{AV}$ .

## 3. Testing and adjusting the approach

Testing was performed using the schemes described in the previous section: HR-minmod, HR-vanLeer, HR-MC, RK2-MC, and RK3-WENO5. The Riemann problem was generally solved by the exact (Godunov) solver [14]; sample computations with the Roe [15] and HLLC [16] solvers gave close results (like in the case of the first-order schemes [3]). In Appendix A we provide a convergence study to demonstrate the efficiency of the selected schemes in the case of solving a shock-free problem (the artificial viscosity approach has no impact on the solution in that case).

The vast majority of the test computations were carried out using various modifications of the Quirk problem. The results differed substantially, which complicated their analysis. To make the material of this section easier to comprehend, let us communicate two important findings of the work in advance: (1) of all the schemes selected, HR-minmod demonstrates the lowest disposition to the generation of spurious oscillations behind the shock; (2) provided that two recommendations to the schemes (to be formulated in this section) are followed, the value of the coefficient $C_{AV} = 0.5$ in the artificial viscosity model ensures suppression of the carbuncle instability for all the schemes under consideration.

### 3.1. Double Mach reflection problem: characteristic variables vs. primitive variables

Before moving on to the parametric testing based on the Quirk-type problems, let us consider the results of computations of one well known test problem, namely the double Mach reflection problem [17]. Here we use the original formulation of the problem with one improvement: at each time step the exact values of flow variables behind the Mach 10 shock are placed at the grid points above and to the left of all shock fronts (specifically, located at a distance of at least five grid intervals from them). The purpose of this is to suppress small spurious perturbations, resulting from inappropriate (inconsistent with the numerical viscosity) smearing of the main shock in the initial data and at the upper boundary of the domain (see section 15.8.4 of the book [6] about the "start-up errors").

Fig. 1 shows the results of computations by the HR-MC scheme on a square grid with spacing $h$ = 1/480. Four variants of the computation that differ in the type of variables used in the reconstruction and presence/absence of artificial viscosity are presented. Fig. 1(a) corresponds to a computation with primitive variables without adding the artificial viscosity. The arrows in the figure point at the perturbations of different forms that noticeably debase the data obtained. The perturbations indicated by arrows 1 and 2 are found behind the reflected shock wave to the left of the triple point. We can say with confidence that they are not caused by the carbuncle instability, because the flow lines in this part of the flow are not aligned with the grid. Other perturbations are found in the Mach stem (the shock between the triple point and the plate $y = 0$). Here we see a sawtooth distortion of the Mach stem (arrow 3) and its kinking (arrow 4) resulting in the non-physical flow structure behind the shock. These perturbations are the manifestations of the carbuncle instability.

Fig. 1(b) shows data obtained using the primitive variables, but with adding the artificial viscosity. We see that this variant of the computation is completely free of the carbuncle manifestations. As for the perturbations of the first and second kind (arrows 1 and 2), they are still present. Although the former become significantly weaker, the latter are as prominent as before. Fig. 1(c) corresponds to the modification with using the characteristic variables and without adding the artificial viscosity. Its comparison with Fig. 1(a) indicates that discarding the primitive variables in favor of the characteristic ones makes the second kind of perturbations completely vanish and noticeably fades out the first kind of perturbations, while the carbuncle manifestations are not affected. Finally, data of the computation with the characteristic variables and artificial viscosity (Fig. 1(d)) demonstrate the best quality: all the spurious perturbations are either removed completely or reduced to a negligibly low level.

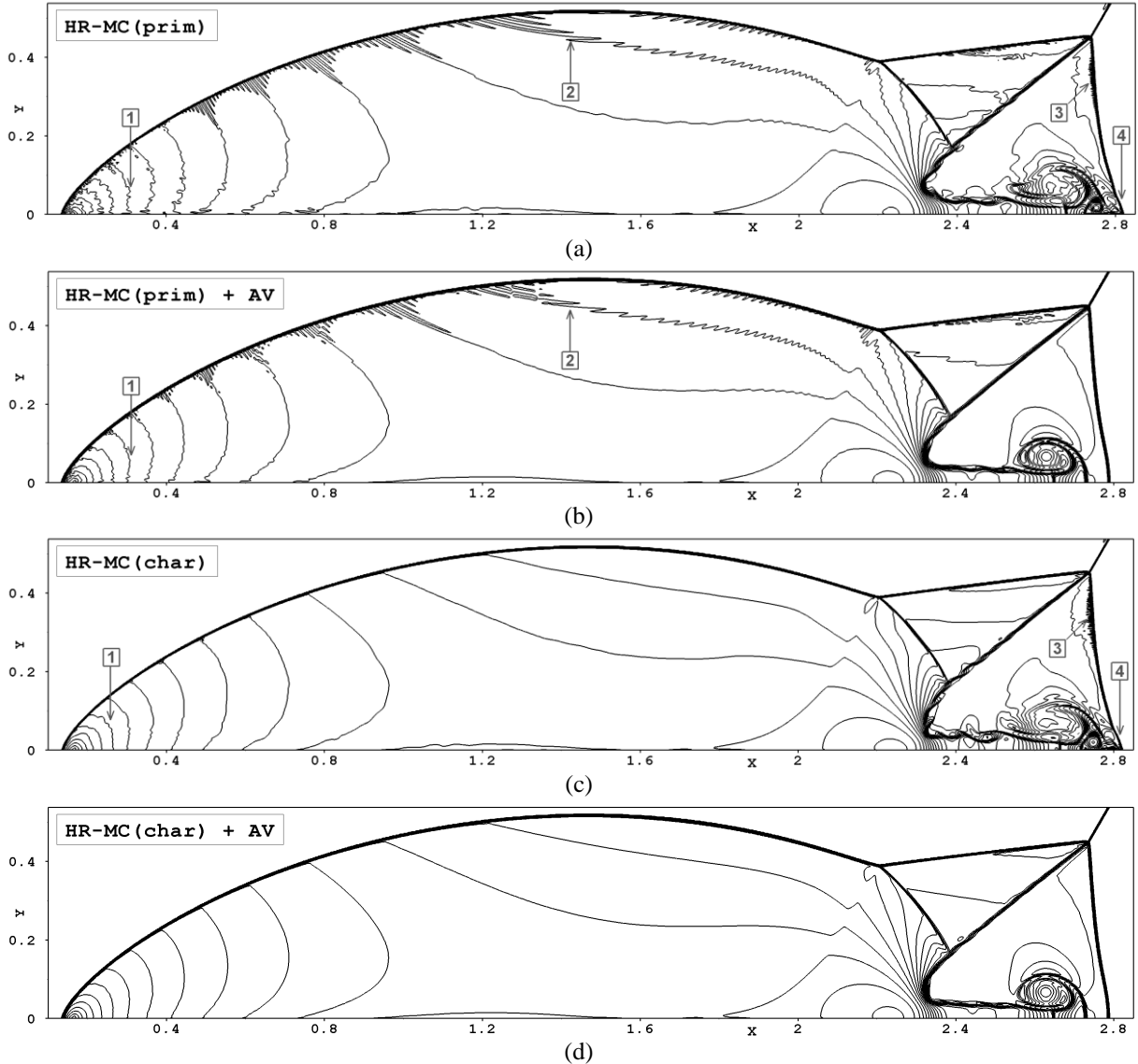

Fig. 1. Double Mach reflection problem. Forty one density contours equally spaced from 2.1 to 22. Computations by the HR-MC scheme on the grid with $h = 1/480$: (a) and (b) reconstruction based on primitive variables; (c) and (d) reconstruction based on characteristic variables; (b) and (d) computations with artificial viscosity.

Thus, we conclude that the artificial viscosity as applied to high-order schemes enables efficient removal of the carbuncle-related spurious perturbations. However, adding it to suppress a different type of perturbations can turn out to be insufficient or ineffective at all. In this case, the problem must be solved by improving the computational techniques themselves. One suggestion for such an improvement can be formulated based on the data presented above; it is not the finding of this work, but is a known recommendation (for example, it can be found in the review article by van Leer [5]).

*Suggestion 1 (using the characteristic variables).* In high-order schemes, it is preferable to perform the procedure of data reconstruction using the characteristic variables.

Thus, all the computations reported below will include Suggestion 1, i.e., in each of the schemes under consideration, the data reconstruction will be performed using the characteristic variables. The point to emphasize is that the use of the primitive variables in the minmod limiter generally produces no visible perturbations behind the shock wave, as it is usually the case for other, less dissipative, reconstructions; nevertheless, Suggestion 1 will also be taken into account for this reconstruction, unless otherwise stated.

### 3.2. Formulation of the Quirk-type test problems

We now consider the Quirk-type test problems to check the stability of a planar shock wave during its computation in the two-dimensional formulation. In [3], three types of shock waves were distinguished in the first-order scheme studies: advancing, reflected and steady shock wave. In this work we consider only the advancing- and steady-shock-wave tests, because they were found to be the most critical.

*Advancing shock wave.* The computational domain of $[0, X] \times [0, Y]$ in the $xy$ plane is covered with a regular grid composed of square cells ($h_y = h_x = 1$) or rectangular cells ($h_y \neq h_x = 1$). Initial conditions of the gas (velocity components, density and pressure) in the domain are the following: $(u_x, u_y, \rho, p) = (0, 0, 1, 1)$. The upper and bottom boundary conditions (at $y = 0$ and $y = Y$) are periodic. The right boundary ($x = X$) is the solid wall, and the boundary condition on the left side ($x = 0$) is the inflow with the values of ($u_1$, 0, $\rho_1$, $p_1$) determined by the shock-wave Mach number ($M_S$) and the ratio of specific heats ($\gamma$):

$$u_1 = u_S \frac{2\left(M_S^2 - 1\right)}{\left(\gamma + 1\right) M_S^2}, \qquad \rho_1 = \frac{\left(\gamma + 1\right) M_S^2}{\left(\gamma - 1\right) M_S^2 + 2}, \qquad p_1 = \frac{2\gamma M_S^2 - \left(\gamma - 1\right)}{\left(\gamma + 1\right)},$$

where $u_S = \sqrt{\gamma} M_S$ is the shock speed.

The instability of the advancing shock wave is triggered by a negligibly small perturbation of the grid line in the cross section $x = 10$: $\hat{x}_{i0, j} = x_{i0, j} + \delta\left(2 \cdot RND_j - 1\right)$, where $i$ and $j$ are the grid indices in the longitudinal and transverse directions, respectively ($0 \leq i \leq I$ and $0 \leq j \leq J$), $i0 = 10$, $\delta = 10^{-4}$, $RND_j$ are random numbers generated in the range of [0, 1]. Thus, the shock wave acquires a single (at the time it crosses the grid line $i = i0$) perturbation over the entire shock layer. The shock instability is revealed by the growth of the following measure of solution deviation from one-dimensional flow:

$$\varepsilon_0 = \max_{i, j}\left(\left| \rho_{i, j} - \overline{\rho}_i \right|\right), \qquad \overline{\rho}_i = \frac{1}{J} \sum_{j=1}^{J} \rho_{i, j}, \tag{6}$$

*Steady shock wave.* In this test problem, the shock front lies on the grid line $i = i0 = I - 10$ (in the vicinity of the right boundary), which corresponds to the cross section $x = x_0 = X - 10$. Accordingly, initial data are the following:

$$\left(u_x, u_y, \rho, p\right) = \begin{cases} \left((u_1 - u_S), \ 0, \ \rho_1, \ p_1\right), & \text{if } x < x_0, \\ \left(-u_S, \ 0, \ 1, \ 1\right), & \text{otherwise.} \end{cases}$$

During the whole computation we will keep unchanged the flow values at the right and left boundaries of the domain. The upper and bottom boundary conditions are periodic, like in the previous problem. The shock instability is also triggered by a grid perturbation, but this time a few grid lines located near the shock front (where $|i - i0| \leq 2$) are perturbed, and $\delta = 10^{-6}$.

*Computations on parallelogram grids.* A number of test computations will be done on parallelogram grids. In these cases, the vertical grid lines $x = \text{const}$ are replaced with slanting lines $x + \alpha y = \text{const}$, and the free parameter $\alpha$ is selected in such a way that the quantity $N_{shift} \equiv \alpha J\left(h_y / h_x\right)$ is a positive integer. The upper boundary is then shifted to the left relative to the lower boundary by $N_{shift}$ cells, and the boundary conditions are still periodic, but with allowance made for the shifting. (Fig. 9 in [3] shows a fragment of such a grid for $\alpha = 1$, $h_y / h_x = 1/2$.) For the advancing shock wave computation, in the cells, the centers of which are located in the region $x < 0$, the inflow values ($u_1$, 0, $\rho_1$, $p_1$) will be set and kept unchanged during the whole computation.

In all the test cases (unless otherwise specified) we assume that $\gamma = 1.4$, and $M_S = 20$ (the most prominent shock instability is observed in hypersonic flows; see [3]). Computations by the RK3-WENO5 scheme will be done with the Courant number $C_{cfl} = 0.6$; for all the other schemes we take $C_{cfl} = 0.8$.

### 3.3. Advancing-shock-wave test problem: rectangular grids

This test problem was computed using different schemes for the following cell aspect ratios: $h_y / h_x = 1/8, 1/4, 1/2, 1, 2, 4, 8$. The number of cells in the domain was $I{\times}J = 1600{\times}50$ for $h_y / h_x \geq 1/2$ and $I{\times}J = 1600{\times}100$ for $h_y / h_x < 1/2$. Each test case was computed repeatedly: we varied the coefficient $C_{AV}$ in the artificial viscosity model and found the minimum value $C_{AV}^{\min}$ to suppress the carbuncle instability. The instability was identified based on the following criterion (below we call it *the* $\varepsilon_0$ *-criterion*): if the value of $\varepsilon_0$ (formula (6)) remains negligible (at the background level), the computation is stable, while its exponential growth is indicative of the instability.

The computations show that in the case of using the HR-minmod and HR-vanLeer schemes, the value of $C_{AV} = 0.20$ ensures the suppression of carbuncle instability in the whole range of the cell aspect ratios. As to the HR-MC and RK3-WENO5 schemes, this statement is valid only for $h_y / h_x \geq 1$, and as to the RK2-MC scheme, for $h_y / h_x \geq 4$. In the other cases, applying the $\varepsilon_0$ -criterion to identify the value of $C_{AV}^{\min}$ was found to be inefficient.

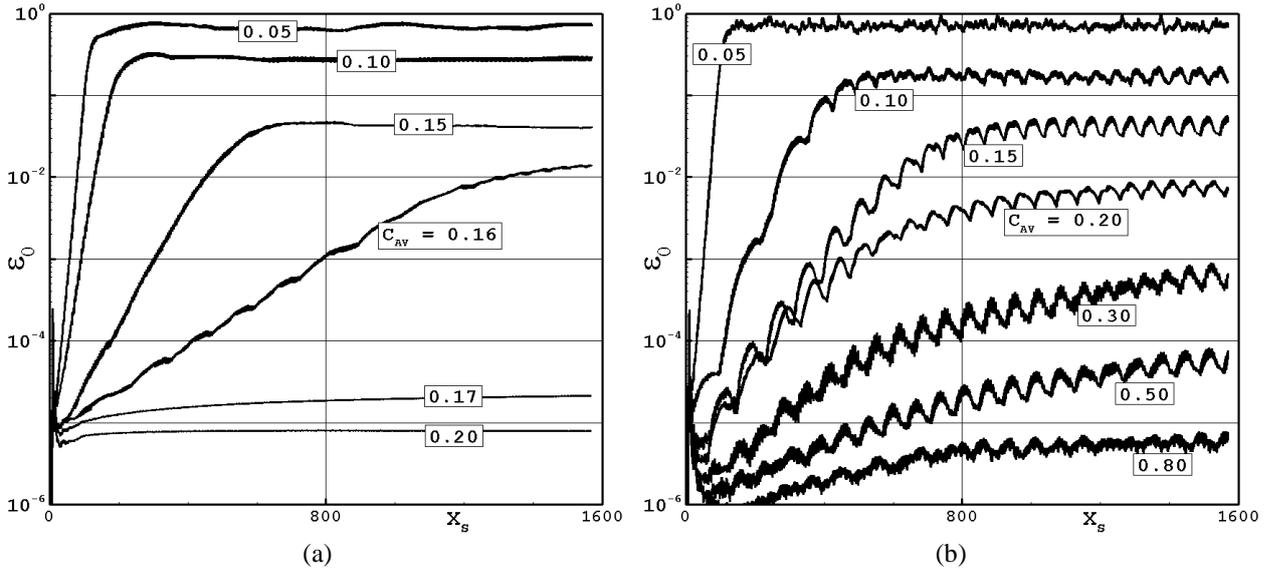

(a)             (b)

Fig. 2. Instability evolution by the $\varepsilon_0$ -criterion in the modified Quirk problem (advancing shock wave; $M_S = 20$). Computations using the HR-MC scheme: (a) $h_y / h_x = 1$; (b) $h_y / h_x = 1/4$.

Fig. 2 shows the plots of $\varepsilon_0$ as a function of distance traveled by the shock wave ($x_S = u_S \cdot t$) for the HR-MC scheme at two values of $h_y / h_x$. In Fig. 2(a) (the case with $h_y / h_x = 1$) one can see that the shock instability attenuates quickly with increase in $C_{AV}$ until completely suppressed; here one can assume that $C_{AV}^{\min} \approx 0.17$. For $h_y / h_x = 1/4$ (Fig. 2(b)), the slope of the curves goes down very slowly with increasing $C_{AV}$, and the curves themselves are highly oscillating. The $\varepsilon_0$ -criterion appears to be inefficient here, hence it is necessary to develop an alternative approach for determining $C_{AV}^{\min}$. This problem was resolved in the following way.

One of the signs of the instability is the presence of overshoots in density behind the shock front. With increasing $C_{AV}$, such overshoots fade out until small enough to be neglected. In this case, one can ignore the residual instability even if it is present in the computation.

Thus, the degree of solution instability behind the shock wave can be quantified as

$$\varepsilon_1 = \max_{i,j}\left(\varepsilon_{i,j}\right), \qquad \varepsilon_{i,j} = \frac{\rho_{i,j}}{\rho_1} - 1,$$

where $\rho_1$ is the exact value of density behind the shock wave, $\varepsilon_{i,j}$ is the local density error. Suppose that the shock instability can be ignored if the value of $\varepsilon_1$ during the computation does not exceed the threshold level of $10^{-3}$. We call this condition *the $\varepsilon_1$-criterion*.

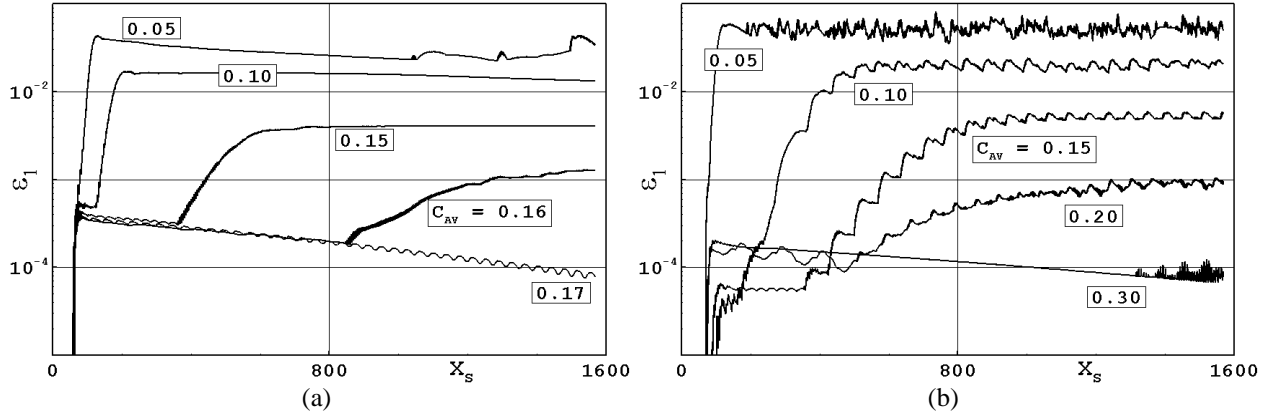

Fig. 3. Instability evolution by the $\varepsilon_1$-criterion in the modified Quirk problem (advancing shock wave; $M_S = 20$). Computations using the HR-MC scheme: (a) $h_y / h_x = 1$; (b) $h_y / h_x = 1/4$.

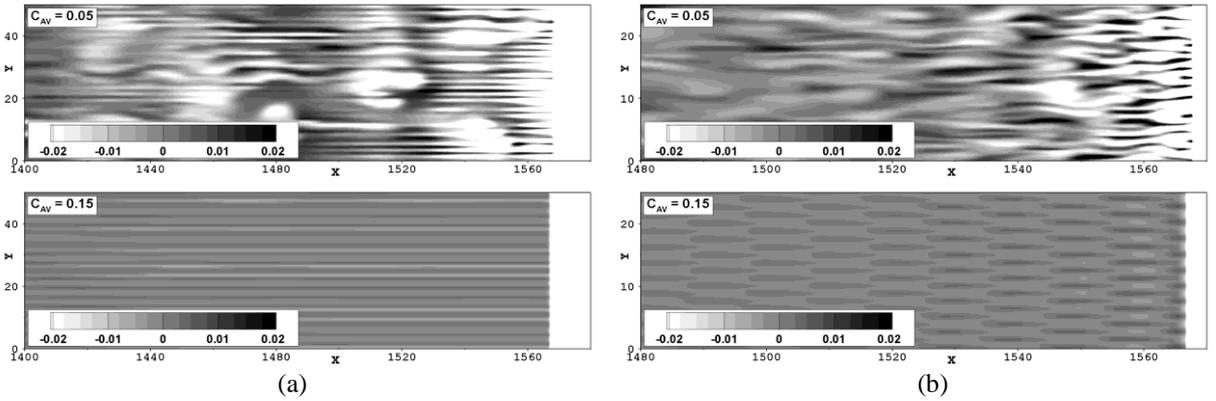

Fig. 4. Distributions of the local density error ($\varepsilon$) in the modified Quirk problem (advancing shock wave; $M_S = 20$). Computations using the HR-MC scheme: (a) $h_y / h_x = 1$; (b) $h_y / h_x = 1/4$.

Fig. 3 shows the plots of $\varepsilon_1(x_S)$ for the same test cases as those chosen for Fig. 2. One can see that the $\varepsilon_1$-criterion on the grid with $h_y / h_x = 1$ (Fig. 3(a)) gives the value of $C_{AV}^{\min}$, which is close to the value obtained with the $\varepsilon_0$-criterion: $C_{AV}^{\min} \approx 0.17$. For $h_y / h_x = 1/4$ (Fig. 3(b)), the new criterion gives $C_{AV}^{\min} \approx 0.20$, which is just a little higher.

It is also of interest to consider the spatial distribution of local errors $\varepsilon_{i,j}$: Fig. 4 shows data characterizing the flow instability in the selected test cases for two values of $C_{AV}$. One can see that at $C_{AV} = 0.05$, the instability behind the shock wave is irregular, and the density error exceeds 2% in both test cases. At $C_{AV} = 0.15$, the instability attenuates to ~0.5% and becomes regular. In the test case with $h_y / h_x = 1$, the shock front has a sawtooth structure (with a period of ~$3h_y$), which does not evolve as the computation progresses (local errors are preserved along the flow lines). In the case with $h_y / h_x = 1/4$, the shock front has a sinusoidal structure (with a period of ~$10h_y$), and the cellular pattern of the downstream density levels points at the cyclic time evolution of the shock wave structure itself rather than at acoustic waves (pressure errors here are an order smaller than density errors).

Thus, the use of the $\varepsilon_1$-criterion in this test problem turned out to be quite effective. It allowed us to establish that the value of $C_{AV} = 0.30$ enables the suppression of the carbuncle instability for any of the schemes under consideration in the entire range of cell aspect ratios: $1/8 \leq h_y / h_x \leq 8$.

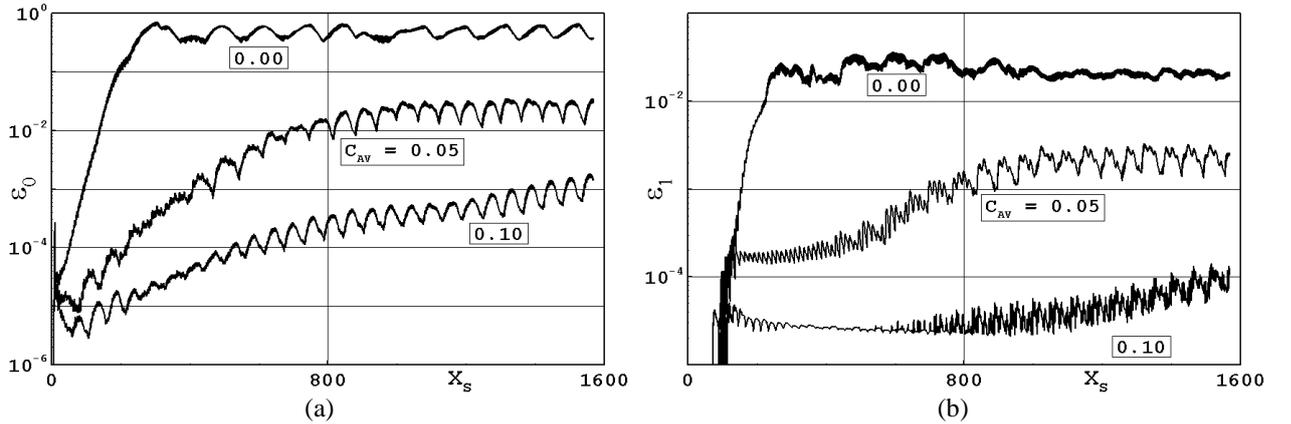

Fig. 5. Instability evolution in the modified Quirk problem (advancing shock wave; $M_s = 20$). Computations using the HR-MC scheme with the HLL solver on the grid with $h_y / h_x = 1/4$: (a) $\varepsilon_0$-criterion; (b) $\varepsilon_1$-criterion.

*On carbuncle-free Riemann solvers.* It is thought that the shock-capturing methods that employ the incomplete Riemann solvers (based on coarse wave models that do not capture the contact discontinuity) are free of the carbuncle instability. However, this is not exactly true. For example, advancing-shock-wave test computations by the HR-MC scheme with the HLL solver [18], which is incomplete, exhibit the carbuncle instability at $h_y / h_x \leq 1/4$. Fig. 5 shows data obtained on the grid with $h_y / h_x = 1/4$. One can see that in the case of $C_{AV} = 0$ the instability develops rather actively, whereas adding the artificial viscosity causes its rapid decay. The appearance of the instability in this case can be explained as follows.

It was established that the cause of the carbuncle instability lies in the insufficient dissipation via the contact discontinuity in the direction parallel to the shock. This happens, if a complete Riemann solver is used on the cell faces aligned with the flow lines. If a high-order scheme is used on a grid with $h_y / h_x \ll 1$, the lateral dissipation decreases substantially promoting the instability even if an incomplete solver is applyed. To this we can add that the computation using the HLL first-order scheme on a grid with $h_y / h_x = 1/32$ (a fairly exotic test case) also displayed a shock instability.

### 3.4. Advancing-shock-wave test problem: parallelogram grids

In this test problem, the advancing shock wave is computed on a parallelogram grid with $\alpha = 1$ (vertical grid lines are replaced with slanting lines $x + y = $ const) and $h_y / h_x = 1, 1/2, 1/4, 1/8$; the grid is composed of $I \times J = 1200 \times 48$ cells. As noted in [3], the numerical solution of the problem remains quasi-one-dimensional in the sense that the solution values at the cell centers located along the same vertical line ($x = $ const) are identical. The shock instability, however, still shows up in this test problem: introducing a small grid perturbation with time leads to significant distortions in the shock layer and downstream. Adding the artificial viscosity in this case also suppresses the instability.

Table 1 shows the results of processing the test problem computations. It lists the minimum values of $C_{AV}$ to suppress the instability based on the $\varepsilon_1$-criterion for the different schemes. One can see that, as applied to the HR-minmod scheme, the value of $C_{AV} = 0.15$ ensures suppression of the instability in the whole range of the cell aspect ratios. (Also note that the use of the $\varepsilon_0$-criterion for this scheme gives close values of $C_{AV}^{\min}$.) The HR-vanLeer and RK3-WENO5 schemes on the whole demonstrate acceptable results, while the HR-MC and RK2-MC schemes encounter the problem on a grid with $h_y / h_x = 1$.

Table 1. The values of $C_{AV}^{\min}$ in the advancing-shock-wave test problem on the parallelogram grids ( $\varepsilon_1$-criterion, $M_S = 20$). The bracketed values correspond to the case of using Suggestion 2 (switching to the minmod limiter).

| $h_y / h_x =$ | 1 | 1/2 | 1/4 | 1/8 |
|---|---|---|---|---|
| HR-minmod | 0.05 | 0.10 | 0.15 | 0.15 |
| HR-vanLeer | 0.20  (0.05) | 0.15  (0.11) | 0.22  (0.19) | 0.27  (0.17) |
| HR-MC | > 1  (0.05) | 0.13  (0.12) | 0.28  (0.25) | 0.45  (0.38) |
| RK2-MC | > 1  (0.07) | 0.30  (0.22) | 0.37  (0.30) | 0.47  (0.43) |
| RK3-WENO5 | 0.40  (0.12) | 0.15  (0.15) | 0.25  (0.25) | 0.25  (0.25) |

Fig. 6(a) shows the plots of $\varepsilon_1(x_S)$ in the HR-MC computations on a grid with $h_y / h_x = 1$. One can see that the level of shock instability decreases very slowly with increasing $C_{AV}$, so the stability criterion $\varepsilon_1 < 10^{-3}$ turns out to be unattainable at moderate values of $C_{AV}$. Consequently, adding the artificial viscosity in this test case appears to be insufficient for suppressing the shock wave instability, and one needs to improve the computational method itself. A technique that was found to be rather efficient for this purpose is based on the following reasons.

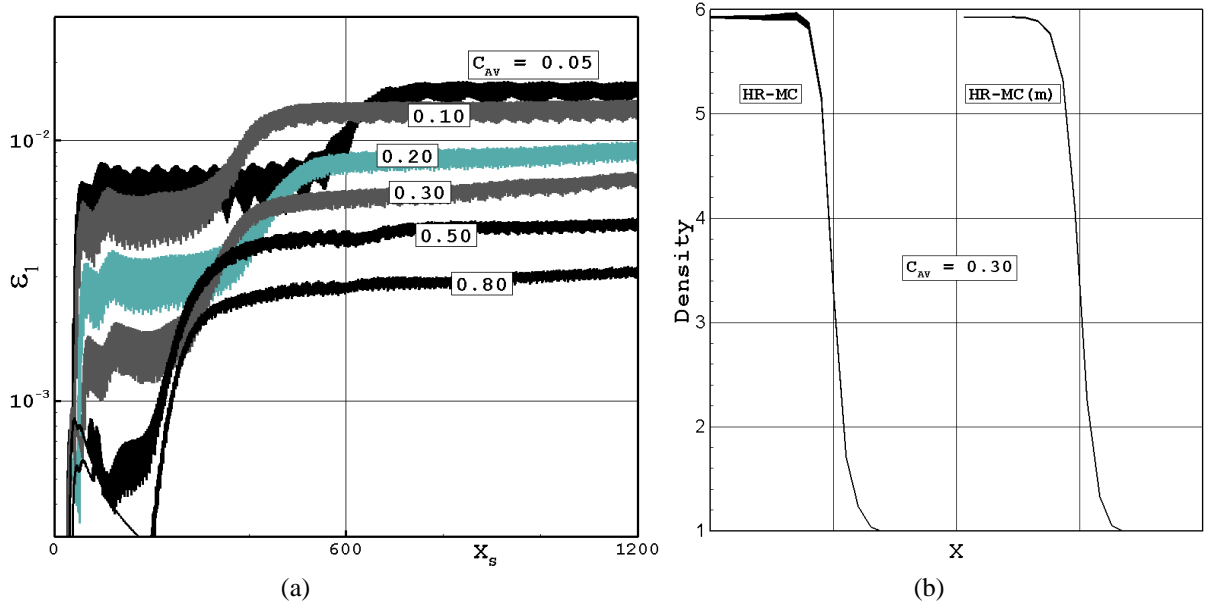

(a)                                                                (b)

Fig. 6. Advancing-shock-wave test computations by the HR-MC scheme on the parallelogram grid ( $h_y / h_x = 1$; $M_S = 20$): (a) instability evolution by the $\varepsilon_1$-criterion ( $C_{AV}$ is varied); (b) shock wave profiles obtained by the basic and modified schemes ( $C_{AV} = 0.3$).

Let us direct our attention to the reconstruction based on the minmod limiter. Although it has the highest dissipation among the reconstructions under consideration, when combined with the artificial viscosity, it suppresses the shock wave instability most efficiently. From the practical viewpoint, the use of more accurate reconstructions is preferable in smooth parts of the solution and at contact discontinuities. However, at the computational points inside the shock layer, the accuracy of reconstruction becomes second to the capability of the computational algorithm to prevent oscillations in the solution. This consideration leads us to the following suggestion for data reconstruction improvement.

*Suggestion 2 (switching to the minmod limiter).* During a reconstruction in the cell $(i, j)$ in one of the grid directions (for example, along the $i$–th direction) we check whether the cell or its neighbors (cells $(i-1, j)$ and $(i+1, j)$) are situated within the shock layer; as an indicator of this event we use the condition $\mu_{AV} \neq 0$. Whenever this condition is fulfilled, we assume that the data distribution in the cell is linear and use the minmod limiter. To denote this technique, we add the character (m) to the basic reconstruction, for example: MC(m).

Fig. 6(b) presents a comparison of the shock wave profiles obtained by the HR-MC scheme (reference case) and the HR-MC(m) scheme (modified scheme) at $C_{AV} = 0.3$. We see that the one-dimensional flow pattern on the top of the shock layer is distorted and a noticeable density overshoot (~0.7%) is present in the reference case. The use of Suggestion 2 (switching to the minmod limiter) results in a significant improvement of the solution: it recovers the one-dimensional flow pattern and suppresses the density overshoots by an order of magnitude. Here, it is important that using the modified scheme does not entail any noticeable shock broadening.

Table 1 in brackets presents the data obtained by the modified schemes. One can note that the use of Suggestion 2 improves the results significantly: the value of the coefficient $C_{AV} = 0.43$ now suppresses the instability in all the cases, and in the absence of highly elongated cells ($h_y / h_x = 1/8$) the value of $C_{AV} = 0.30$ turns out to be sufficient.

### 3.5. Steady-shock-wave test problem: rectangular grids

This test problem was computed for the following cell aspect ratios: $h_y / h_x = 1/4$, 1/2, 1, 2. The domain initially contained $I \times J = 400 \times 100$ cells.

Fig. 7 shows data obtained by the HR-minmod scheme on a grid with $h_y / h_x = 1/2$. It demonstrates the instability evolution with respect to the $\varepsilon_0$-criterion and the spatial distribution of local errors for two values of $C_{AV}$. One can see that the value of the coefficient $C_{AV} = 0.25$ at early times ($t < 60$ in Fig. 7(a)) suppresses the shock-wave instability, but as the computation progresses with this value of $C_{AV}$, the instability still appears with large-scale structures behind the shock front (Fig. 7(b)). Such a delayed manifestation of the instability is attributed to the following.

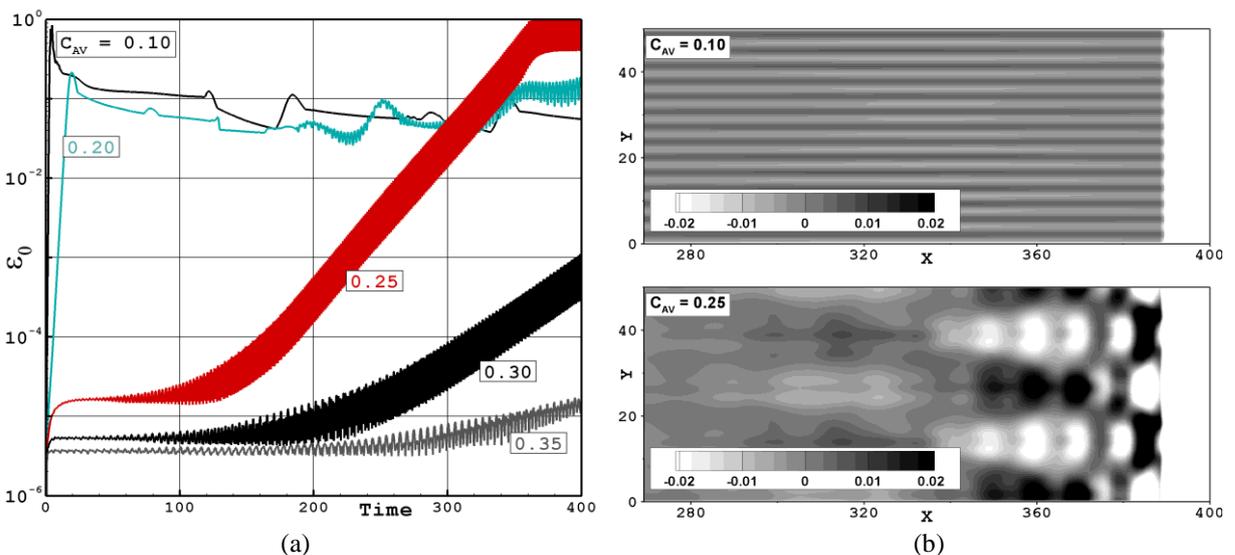

(a)                       (b)

Fig. 7. Steady-shock-wave test computations by the HR-MC scheme ($h_y / h_x = 1/2$; $M_S = 20$): (a) instability evolution by the $\varepsilon_0$-criterion ($C_{AV}$ is varied); (b) distributions of the density error ($\varepsilon$) for two values of $C_{AV}$.

As revealed in [19], the appearance of the carbuncle instability in steady-shock-wave computations may depend on the relative positioning of the shock on the grid. In our case, the initial data correspond to the shock positioning strictly on a grid line $i = i_0$. At the very beginning of the computation, we observe shock broadening due to the artificial viscosity, which causes a flow perturbation that moves downstream. This perturbation reaches the right boundary of the domain, is partly reflected and returns back to the shock wave in time $t = -x_0/(u_x - a) + x_0/(u_x + a) \approx 87$ (here, $x_0 = 390$, $u_x \approx -4$, $a \approx 10.5$). As a result, the shock front shifts a little, which is critical for the appearance of the instability at this value of $C_{AV}$.

To account for the possible sensitivity of the instability to the shock position, this test problem was modified as follows.

*Quasi-steady shock wave*. As distinct from the original problem setup, a small quantity of 0.002 is added to the longitudinal velocity component all over the domain and at its boundaries. In this case, the shock front slowly drifts across the grid traversing one cell in $\Delta t = 500$, which requires 30,000 to 80,000 time steps depending on the particular test case. The grid perturbation increases to $\delta = 10^{-3}$, which makes it easier to reveal the instability once it appears in some phase of a computation (at a certain relative shock position on the grid). The size of the domain is decreased to $I \times J = 100 \times 100$.

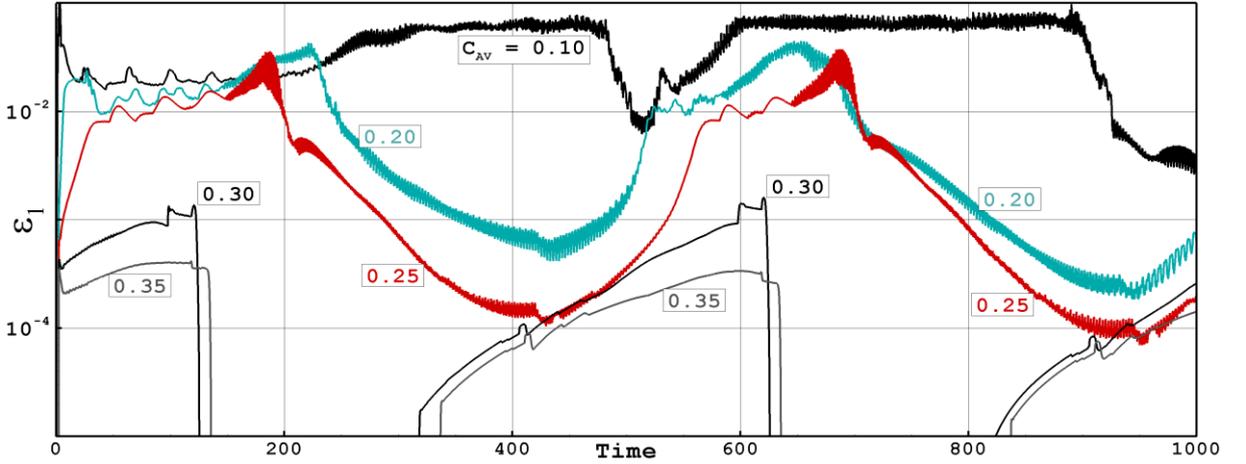

Fig. 8. Instability evolution by the $\varepsilon_1$-criterion in modified Quirk problem (quasisteady shock wave; $M_S = 20$). Computations using the HR-minmod scheme on the grid with $h_y / h_x = 1/2$.

Table 2. The values of $C_{AV}^{\min}$ in the quasisteady-shock-wave test problem on the rectangular grids ($\varepsilon_1$-criterion, $M_S = 20$). The bracketed values correspond to the case of using Suggestion 2 (switching to the minmod limiter).

| $h_y / h_x =$ | 1/4 | 1/2 | 1 | 2 |
|---|---|---|---|---|
| HR-minmod | 0.35 | 0.35 | 0.30 | 0.25 |
| HR-vanLeer | 0.50 (0.35) | 0.50 (0.35) | 0.30 (0.30) | 0.25 (0.25) |
| HR-MC | >1 (0.35) | >1 (0.35) | >1 (0.30) | ~1 (0.25) |
| RK2-MC | >1 (0.35) | >1 (0.35) | >1 (0.35) | ~1 (0.25) |
| RK3-WENO5 | >1 (0.35) | >1 (0.35) | >1 (0.35) | >1 (0.25) |

Fig. 8 shows the instability evolution by the $\varepsilon_1$-criterion in the modified test problem (quasi-steady shock wave) for the same test cases as shown in Fig. 7. It is easy to see that (1) at $C_{AV} > 0.2$, the dependence $\varepsilon_1(t)$ is strictly cyclic with a period of $\Delta t = 500$, and (2) the value of $C_{AV} = 0.35$ satisfies the condition $\varepsilon_1 < 10^{-3}$ at each point of time, i.e., at each shock position relative to the grid. Thus, for the given test case (quasisteady-shock-wave computation by the HR-minmod scheme on a grid with $h_y / h_x = 1/2$), one can assume that $C_{AV}^{\min} \approx 0.35$.

Similarly, the quasisteady-shock-wave test problem was calculated by other schemes and with other cell aspect ratios. We have considered both basic schemes and their modifications using

Suggestion 2 (switching to the minmod limiter). Table 2 shows the results of processing the test problem computations. One can see that it is highly recommended to use Suggestion 2 for this test problem for suppressing the instability (by the $\varepsilon_1$-criterion). The only scheme to demonstrate relative efficiency in its basic form is the HR-vanLeer scheme, while the switching to the minmod limiter makes it possible to suppress the instability at a moderate value of the free coefficient: $C_{AV} = 0.35$.

### 3.6. Steady-shock-wave test problem: parallelogram grids

In this problem, the steady shock wave is computed on parallelogram grids with $N_{shift} = 2$ (the upper boundary of the domain is shifted relative to the lower boundary by two cells) and $h_y / h_x = 1/4$, $1/2$, $1$, $2$; the grid is composed of $I \times J = 200 \times 200$ cells. Here, the grid lines $i = \mathrm{const}$ slightly deviate from the line of the shock front, so the whole spectrum of shock front positions relative to the grid is present in the computation. Note that this type of problems often displays significant distortions behind the shock, which presents a serious computational challenge.

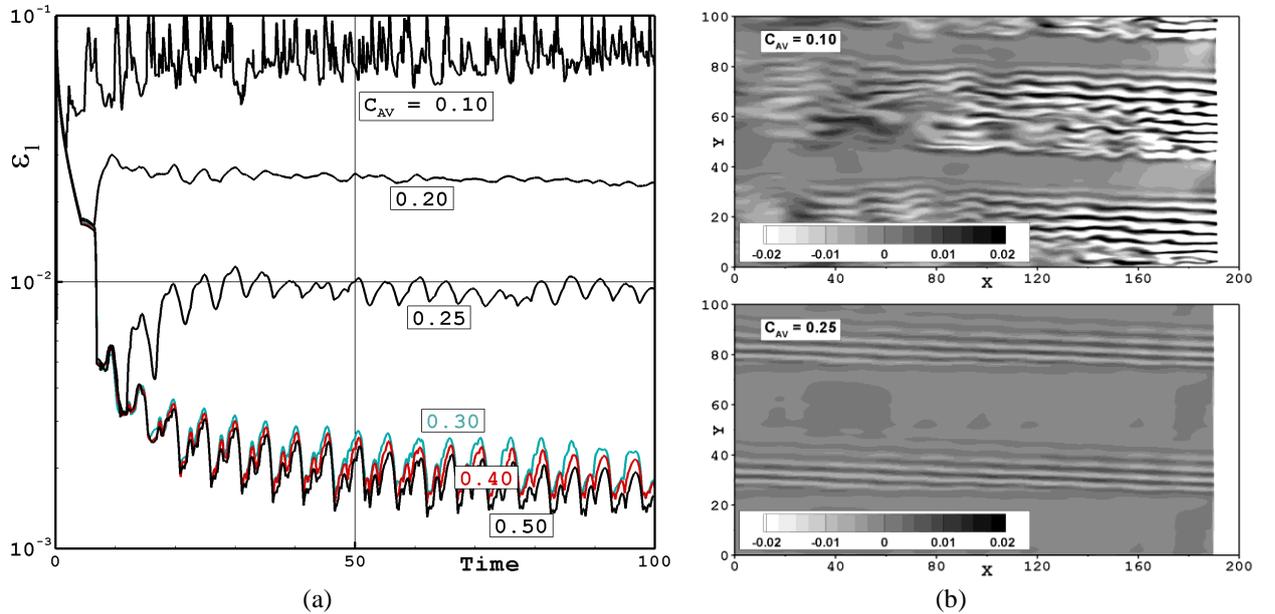

Fig. 9. Steady-shock-wave test computations by the HR-MC scheme on the parallelogram grid ( $h_y / h_x = 1/2$; $M_S = 20$): (a) the time evolution of the maximum density error $\varepsilon_1$; (b) distributions of the density error ( $\varepsilon$ ) at $t = 100$ for two values of $C_{AV}$ .

Fig. 9 shows data obtained by the HR-minmod scheme on a grid with $h_y / h_x = 1/2$. One can see that the errors in density decrease substantially as $C_{AV}$ increases to $C_{AV} \approx 0.3$ (Fig. 9(a)), but any further increase in $C_{AV}$ produces a minor effect. In addition, the following observation was made. At small values of $C_{AV}$, the flow pattern in the upper part of the domain ($50 < y < 100$ in Fig. 9(b)) differs noticeably from the flow pattern in its lower part ($0 < y < 50$), which indicates that the flow is unstable (sensitive to small grid perturbations). At $C_{AV} \geq 0.25$, both upper and lower parts of the flow are nearly identical (the greatest density difference between them does not exceed 0.1%), which evidences that the flow is stable. Thus, $C_{AV} = 0.25$ is sufficient for curing the carbuncle in this test case, but the post-shock overshoots in density are still present above $\varepsilon_1 = 10^{-3}$ at any moderate values of $C_{AV}$.

Similarly, steady-shock-wave computations on a parallelogram grid were conducted using other schemes and cell aspect ratios. The value of $C_{AV} = 0.3$ was found to suppress the carbuncle instability for all schemes under consideration, except RK2-MC, but applying Suggestion 2 (switching to the minmod limiter) to this scheme makes it equally stable.

Table 3 shows the maximum error in density at $C_{AV} = 0.5$ (recommended) for the different schemes and cell aspect ratios; the bracketed values correspond to the modifications with Suggestion 2. One can see that (1) the error in the post-shock density in all the cases does not exceed 1%, and (2) applying Suggestion 2 renders the solution noticeably better (in some cases, the error is reduced by half or more). The table also demonstrates that the flow distortions are much stronger at $h_y / h_x \ll 1$, where the grid cells are highly elongated in the direction normal to the shock. However, at $h_y / h_x \gg 1$ (grid cells are highly elongated in the direction tangential to the shock), there is another trouble: the artificial viscosity that is required for suppressing the carbuncle instability can lead to oversmearing of the shock layer, because in this case $\mu_{AV} \sim h_y^2$. From this one can conclude that it is desirable, where possible, to avoid the use of any grid with highly elongated cells in the vicinity of the shock wave.

Table 3. The maximal values of $\varepsilon_1 \times 10^2$ in the steady-shock-wave test problem on parallelogram grids ($M_S = 20$, $C_{AV}$ = 0.5). The bracketed values correspond to the case of using Suggestion 2 (switching to the minmod limiter).

| $h_y / h_x =$ | 1/4 | 1/2 | 1 | 2 |
|---|---|---|---|---|
| HR-minmod | 0.30 | 0.20 | 0.14 | 0.08 |
| HR-vanLeer | 0.48 (0.36) | 0.28 (0.24) | 0.22 (0.18) | 0.12 (0.10) |
| HR-MC | 0.70 (0.38) | 0.44 (0.26) | 0.28 (0.18) | 0.20 (0.10) |
| RK2-MC | 0.85 (0.40) | 0.58 (0.26) | 0.40 (0.20) | 0.22 (0.10) |
| RK3-WENO5 | 0.82 (0.42) | 0.78 (0.32) | 0.70 (0.22) | 0.46 (0.12) |

### 3.7. Effect of the ratio of specific heats

Let us consider the effect of the ratio of specific heats of the gas. Table 4 presents the results of HR-minmod computations for two typical test cases at different values of $\gamma$. In [3] it has been demonstrated that advancing-shock-wave computations by the Godunov scheme (first-order scheme) reveal strong sensitivity of $C_{AV}^{\min}$ to $\gamma$. Table 4 (Case 1) shows that the dependence of $C_{AV}^{\min}(\gamma)$ relaxes with switching to the second-order scheme, and the values of $C_{AV}^{\min}$ do not exceed 0.25 within the considered range of $1.15 \leq \gamma \leq 2$. For the quasi-steady-state shock wave (Case 2 in Table 4), the dependence of $C_{AV}^{\min}(\gamma)$ is extremely weak: in this case, $C_{AV}^{\min} \approx 0.35$ within the error of 0.05. This leads us to a conclusion that the recommended value of $C_{AV} = 0.5$ can be used in a wide range of values of the ratio of specific heats.

Table 4. The values of $C_{AV}^{\min}$ obtained by the HR-minmod scheme for different values of $\gamma$. Case 1: the advancing-shock-wave test problem with $h_y / h_x = 1$. Case 2: the quasi-steady-shock-wave test problem with $h_y / h_x = 1/2$.

| $\gamma =$ | 1.15 | 1.25 | 1.4 | 1.6667 | 2.0 |
|---|---|---|---|---|---|
| Case 1 | 0.25 | 0.22 | 0.19 | 0.15 | 0.12 |
| Case 2 | 0.35 | 0.35 | 0.35 | 0.35 | 0.35 |

### 3.8. Three-dimensional case

To tune the method under consideration for the three-dimensional case, additional test computations were performed. The computations were carried out for advancing and quasi-steady shock waves on a cubic mesh ($h_x = h_y = h_z$) by the HR-minmod scheme with the artificial viscosity for various values of $C_{AV}$. Comparing the new data with those obtained in the two-dimensional case reveals their close similarity. Thus, for the advancing shock wave the minimum value of $C_{AV}$ to suppress the instability is $C_{AV}^{\min} \approx 0.19$ in both 2D and 3D cases. The new value of $C_{AV}^{\min}$ in the quasi-steady-shock computation was also close to the previous value obtained for the two-dimensional case: $C_{AV}^{\min} \approx 0.30$. Hence we conclude that the recommended value of $C_{AV} = 0.5$ as applied to the high-order schemes can be used for simulating problems of any dimensionality.

## 4. Numerical examples

In this section we apply the selected schemes (HR-minmod, HR-vanLeer, HR-MC, RK2-MC and RK3-WENO5) to several test problems to demonstrate the efficiency of the proposed approach for suppressing the carbuncle phenomenon and post-shock oscillations. In conjunction with the selected schemes we employ various Riemann solvers (RS) such as the exact (Godunov), Roe and HLLC solvers. In what follows it is assumed that (1) all reconstruction procedures are based on the characteristic variables (Suggestion 1), unless otherwise stated, and (2) the artificial viscosity model is applied as described above, namely with the principal coefficient $C_{AV} = 0.5$ and with switching to the minmod limiter in the vicinity of the shock (Suggestion 2). Usually we used a Courant number of 0.6 for RK3-WENO5 and 0.8 for the other schemes. However when simulating the steady-state problems without adding the artificial viscosity, to decrease the magnitude of final unsteadiness, we reduced the normal value of $C_{cfl}$ by a factor two or (in case of using the Roe solver) four.

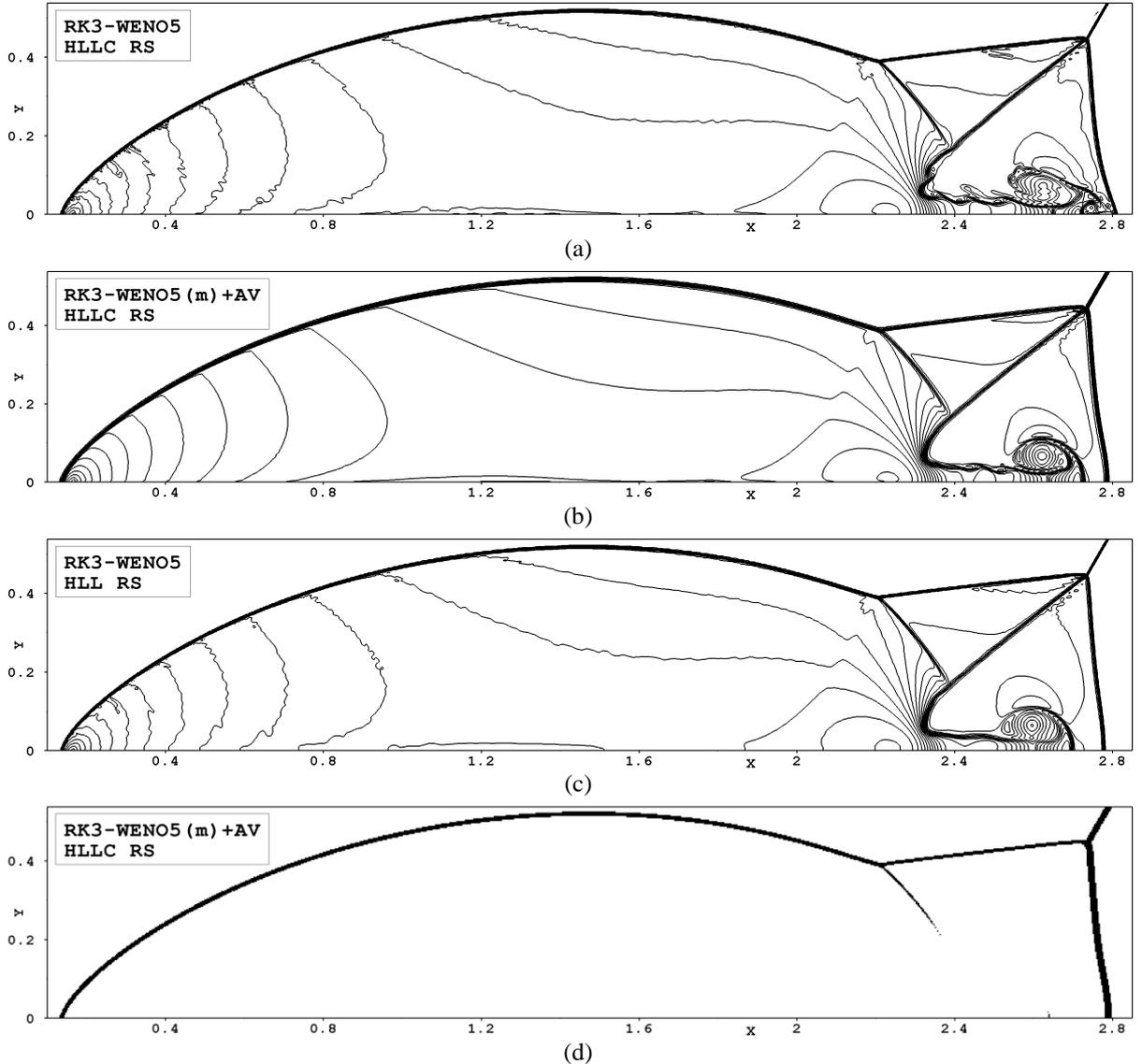

Fig. 10. Double Mach reflection problem. Forty one density contours equally spaced from 2.1 to 22 (a, b, c) and subset of cells (flooded with black color) in which $\mu_{AV} \neq 0$ (d). Computations by the RK3-WENO5 scheme on the grid with $h$ = 1/330: (a) using the HLLC solver without artificial viscosity; (b, d) using the HLLC solver with artificial viscosity; (c) using the HLL solver without artificial viscosity.

## 4.1. Double Mach reflection problem (additional materials)

Let us return to the double Mach reflection problem. In section 3.1 we presented the data obtained by various versions of the HR-MC scheme. To this we can say that the later-developed version of the scheme, called HR-MC(m)+AV, produces data hardly distinguishable from those shown in Fig. 1(d).

Fig. 10 shows the numerical results obtained by the RK3-WENO5 scheme with the HLL-type solvers and the grid spacing $h = 1/330$. (Note that the RK3-WENO5 scheme with $h = 1/330$ consumes nearly twice as much CPU time as the HR-MC scheme with $h = 1/480$.) As one can see, when the HLLC solver is used without artificial viscosity (Fig. 10(a)), the solution is highly distorted by the carbuncle phenomenon and post-shock oscillations. Adding the artificial viscosity (Fig. 10(b)) remedies the defects of the solution completely (note that the additional dissipation is introduced in shock layers only, as illustrated in Fig. 10(d)). On the other hand, applying the HLL solver provides a carbuncle-free solution without adding the artificial viscosity (Fig. 10(c)), but in this case (1) the post-shock oscillations are at a high level and (2) the contact discontinuity, that emerges from the triple point and rolls up to the vortex, seems to be more smeared (compare with Fig. 10(b)).

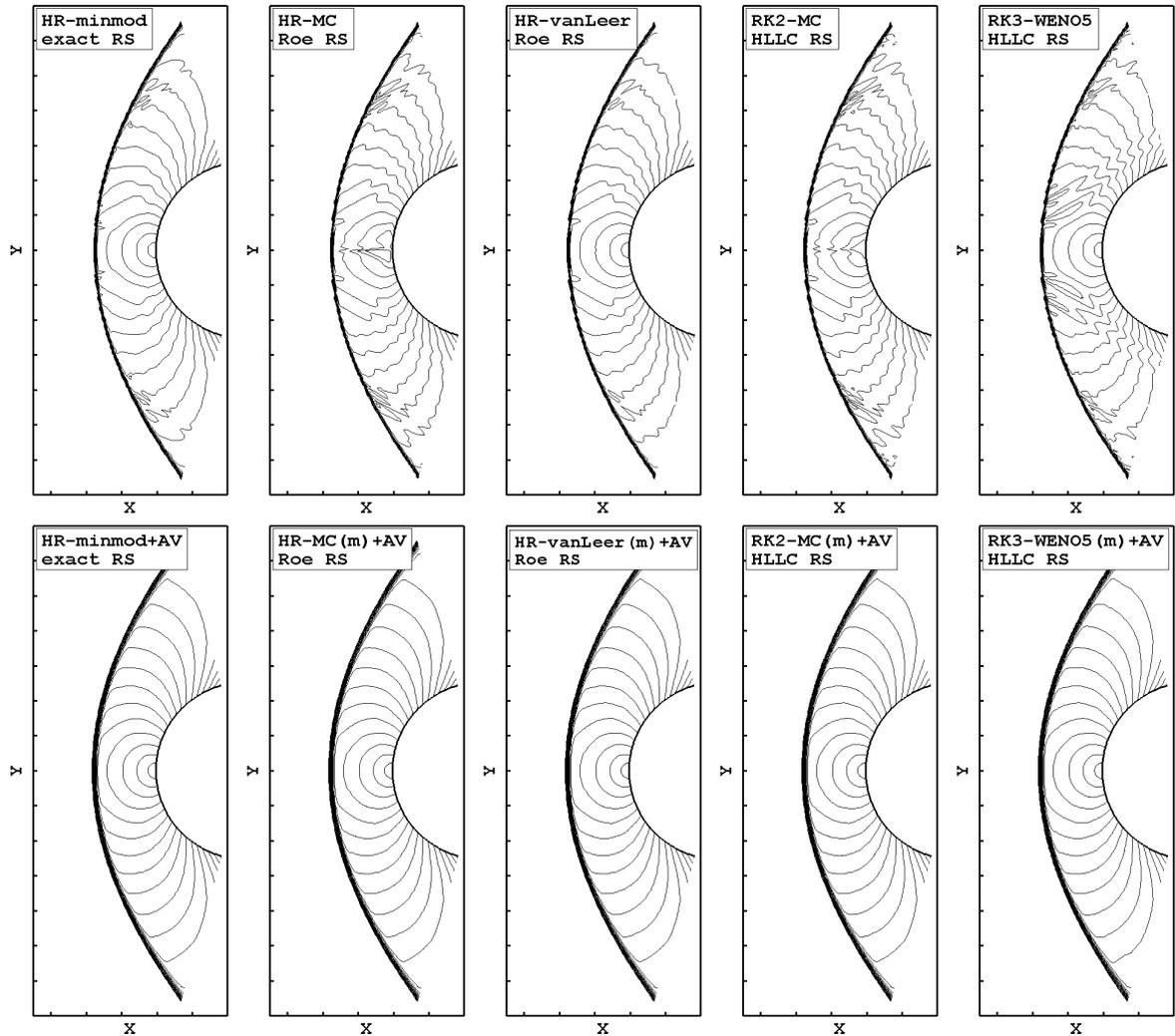

Fig. 11. Supersonic flow around a cylinder at $M_\infty = 3$. Twenty five equidistant Mach number contours from 0.1 to 2.5. Top: original schemes, bottom: respective schemes with artificial viscosity.

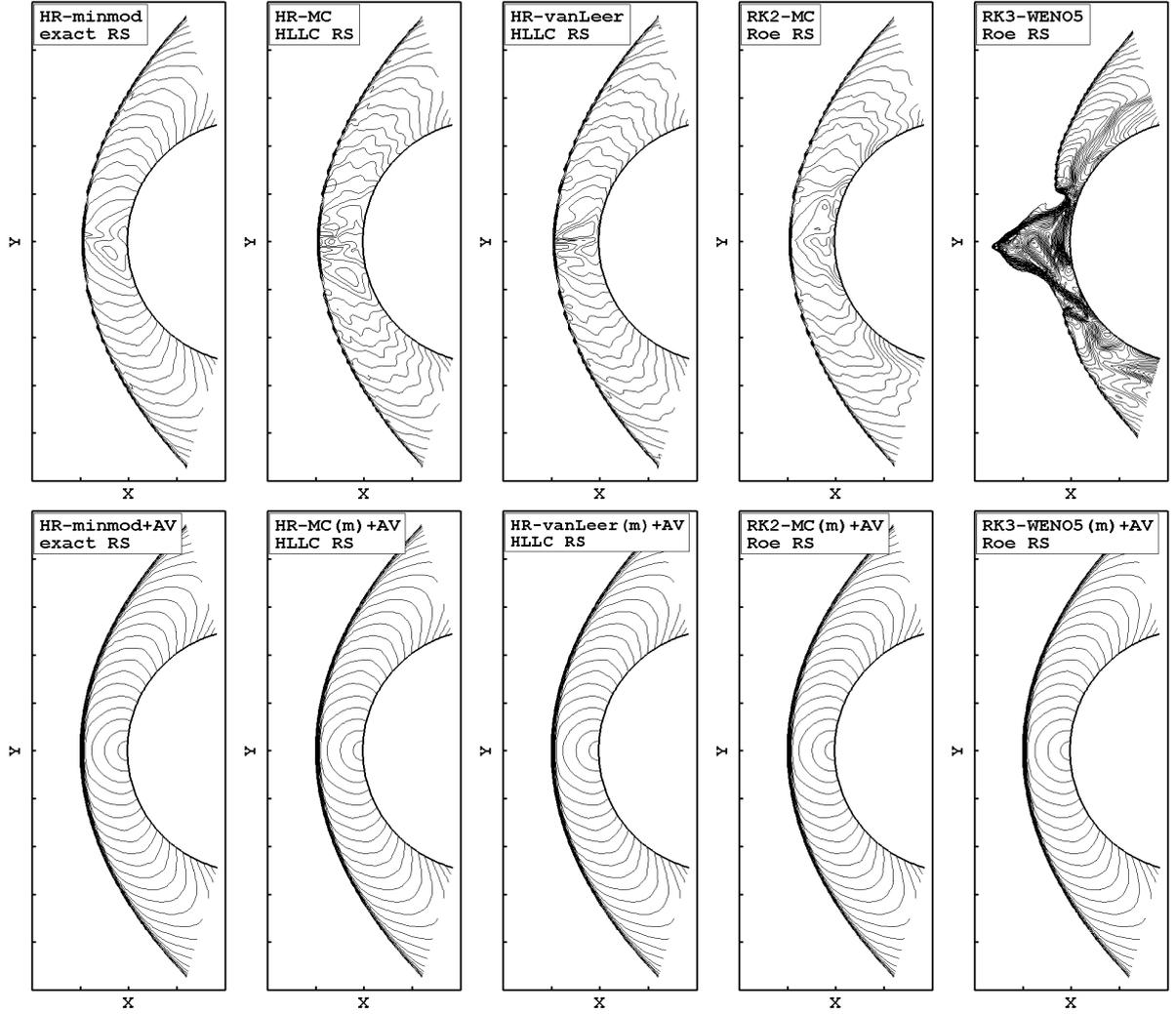

Fig. 12. Hypersonic flow around a cylinder at $M_\infty = 20$. Twenty five equidistant Mach number contours from 0.1 to 2.5. Top: original schemes, bottom: respective schemes with artificial viscosity.

## 4.2. Super- and hypersonic flow around a cylinder (2D case)

For this test problem we use the computational domain $[1, 3]\times[-75°, 75°]$ in polar coordinates $r\varphi$, which is covered with a uniform grid ($J\times I = 90\times160$ cells). The boundary $j = 0$ (the cylinder surface) is a solid wall, and $j = J$ is the inflow with $(u_x, u_y, \rho, p) = (\sqrt{\gamma}M_\infty, 0, 1, 1)$. The bottom ($i = 0$) and upper ($i = I$) boundary condition is a free outflow.

The computed results with $\gamma = 1.4$ and two values of the freestream Mach number, $M_\infty = 3$ and $M_\infty = 20$, are shown in Fig. 11 and Fig. 12, respectively. All the selected schemes in conjunction with various Riemann solvers (RS) were used to simulate these test cases. It can be seen that all the original schemes (i.e. without artificial viscosity) in both test cases suffer from the shock instability (top panels). Adding the artificial viscosity (bottom panels) improves the solutions, so that all the schemes demonstrate similar results.

## 4.3. Heating problem (2D and 3D cases)

Let us now turn to the heating problem that appears when simulating the *viscous* hypersonic flow around a cylinder. To this end we consider the Navier-Stokes equations in both two- and three-dimensional formulations. The flow conditions are the same as those specified in [20]: $M_\infty = 8.1$, the freestream pressure $p_\infty = 370.6$ Pa, the freestream temperature $T_\infty = 63.73$ K, the radius of the cylinder $r = 0.02$ m, the Reynolds number $\text{Re}_r = 1.31 \times 10^5$, the wall temperature $T_w = 300$ K. The working gas is air approximated by the polytropic gas model with $\gamma = 1.4$, the molecular viscosity coefficient is calculated by Sutherland's formula, and the Prandtl number is $\text{Pr} = 0.72$.

In the two-dimensional case we use a computational grid of $J \times I = 150 \times 160$ cells that are uniform in the circumferential direction and varying exponentially in the radial direction: $h_{j+1} / h_j = 1.05$, $h_1 = 5 \times 10^{-5} r$ (as shown in Fig. 13). In the three-dimensional case the grid acquires $K = 45$ uniform intervals in the $z$-coordinate direction; in this direction the domain represents a segment of $0 \le z \le r$ with the periodic boundary condition.

In this test problem we combine the physical (molecular) and artificial viscosities by summing them up. In so doing we prevent generating the artificial viscosity within the boundary layer (which is very thin) by the following condition: if the cell is highly elongated, more specifically if $\max\left(h_i, h_j, h_k\right) > 10 \cdot \min\left(h_i, h_j, h_k\right)$, then $\mu_{AV} = 0$. Besides, we use the local time stepping technique for this problem.

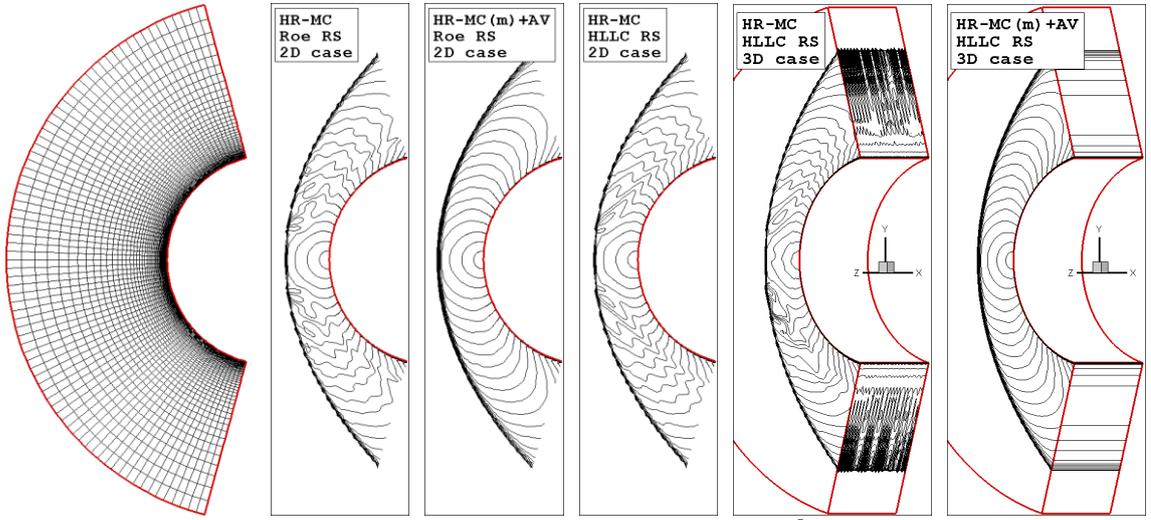

Fig. 13. Viscous hypersonic flow around a cylinder at $M_\infty = 8.1$, $\text{Re}_r = 1.31 \times 10^5$ (2D and 3D cases). Grid (every other grid line) and Mach number contours (from 0.1 to 2.5 with step 0.1).

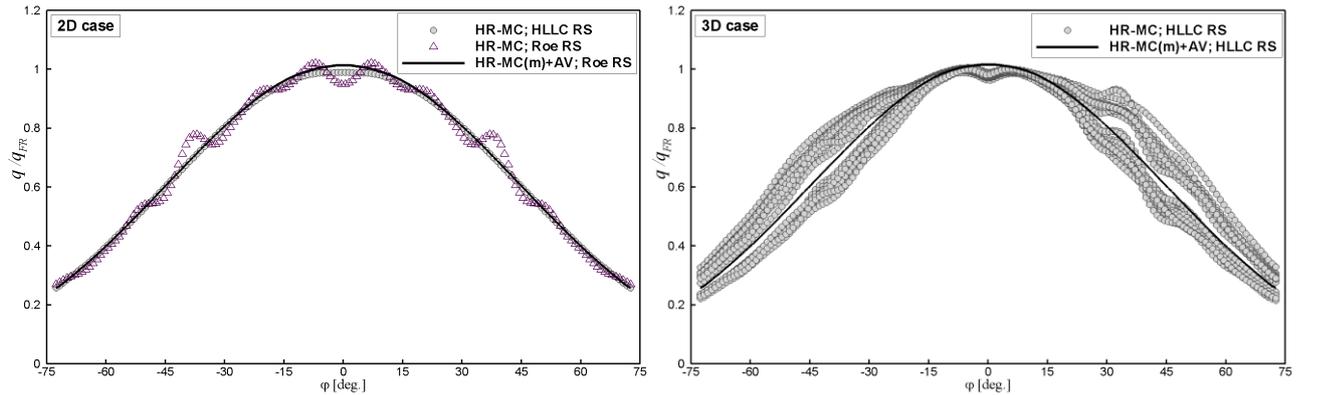

Fig. 14. Heat flux profiles over a cylinder at $M_\infty = 8.1$, $\text{Re}_r = 1.31 \times 10^5$; 2D (left) and 3D (right) cases.

The Mach number contours obtained by the HR-MC scheme with two Riemann solvers (Roe and HLLC) in 2D and 3D cases are shown in Fig. 13. One can see significant flow perturbations in all the computations without artificial viscosity, while applying the artificial viscosity renders the solutions unperturbed. Fig. 14 shows the corresponding wall heating profiles normalized to the Fay–Riddell stagnation value $q_{\text{FR}} = 17.5$ W/cm$^2$. It clearly demonstrates the influence of the post-shock oscillations on the wall heating. Although for the HR-MC scheme with the HLLC solver and without artificial viscosity the profile is not much perturbed in the 2D case, passing to the 3D case deteriorates the data severely. But again, the profiles obtained with the artificial viscosity are of high quality.

## 4.4. Hypersonic flow around a sphere (3D case)

Here we simulate a Mach 10 flow around a sphere using a regular grid composed of $I{\times}J{\times}K = 100{\times}100{\times}36$ cells. The grid point coordinates are set as follows:

$$\mathbf{r}_{ijk} = \left(\frac{1}{2} + \frac{0.24k}{K}\right)\mathbf{n}_{ij}, \quad \mathbf{n}_{ij} = \frac{\mathbf{s}_{ij}}{|\mathbf{s}_{ij}|}, \quad \mathbf{s}_{ij} = \left(-1, \tan\alpha_j, \tan\alpha_i\right), \quad \alpha_j = \frac{5\pi}{12}\left(\frac{2j}{J} - 1\right), \quad \alpha_i = \frac{5\pi}{12}\left(\frac{2i}{I} - 1\right).$$

The boundary $k = 0$ (the sphere surface) is a solid wall, and $k = K$ is the inflow with $(u_x, u_y, u_z, \rho, p) = (\sqrt{\gamma}M_\infty, 0, 0, 1, 1)$. The boundary condition at $i = 0$, $i = I$, $j = 0$ and $j = J$ is a free outflow.

Fig. 15 shows the computational results obtained by the HR-MC scheme with the Roe and HLLC solvers. It can be seen that the simulations without artificial viscosity (top row) demonstrate significant distortion of the solution, especially near the symmetry planes of the grid ($y = 0$ and $z = 0$). Adding the artificial viscosity (bottom row) remedies the defects, and the new solutions become visually similar and cylindrically symmetric.

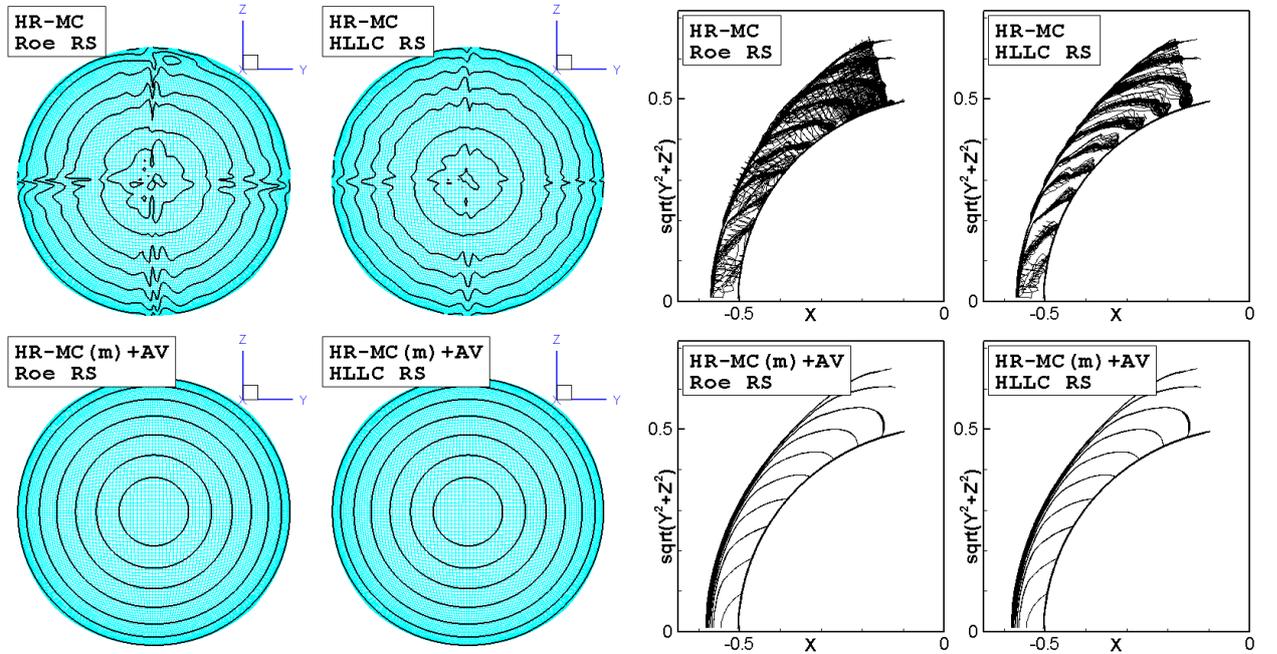

Fig. 15. Results for the Mach 10 flow around a sphere by the HR-MC scheme with the Roe and HLLC solvers. Top row: without artificial viscosity; bottom row: with artificial viscosity. Ten Mach number contours equally spaced from 0.25 to 2.5 on the spherical surface $r = 0.55$ (left; view in $yz$ plane) and on all grid surfaces $j = $ constant (right; view in cylindrical coordinates).

## 4.5. Sedov blast wave problem (2D case)

This test case is a spherically symmetric explosion computed in the 2D formulation. The domain $[0, 1.1]{\times}[0, 1.1]$ in axisymmetric cylindrical coordinates $xr$ is covered with a uniform grid ($J{\times}I = 220{\times}220$). The initial conditions of a polytropic gas with $\gamma = 1.4$ in the domain are $(u_x, u_r, \rho, p) = (0, 0, 1, 10^{-6})$; the only exception is the cell at the origin, where the pressure is $4.334474{\times}10^5$. The bottom and left boundary condition is symmetric flow; and the upper and right boundaries are solid walls. Such a setup corresponds to a total energy of $E_0 = 0.8510719$, at which the shock front reaches a radius of $R_{\text{shock}} = 1$ at time $t = 1$.

Fig. 16 shows the numerical results obtained by the HR-MC and RK3-WENO5 schemes with the Roe solver. As one can see, the outputs of the original schemes (left plots) exhibit carbuncle-like flaws close to both *x* and *r* axes. But then again, adding the artificial viscosity into each scheme (right plots) remedies the problem completely.

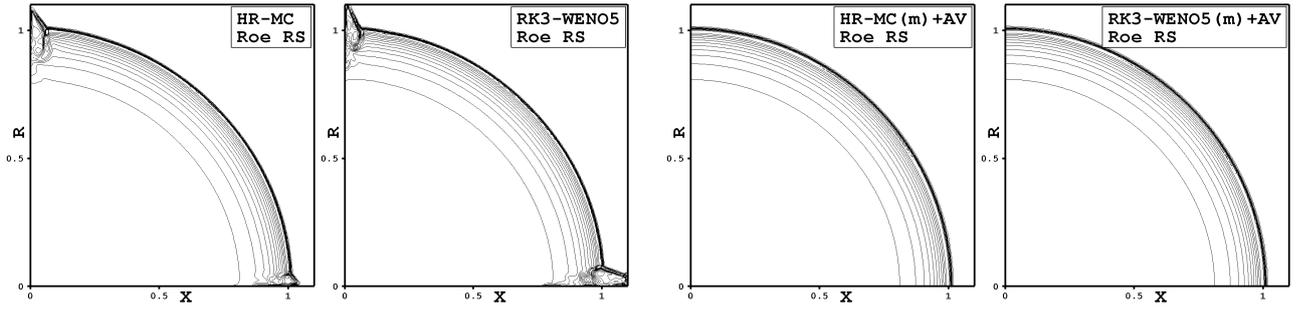

Fig. 16. Sedov blast at *t* = 1. Fifteen density contours equally spaced from 0.4 to 6. Left: original schemes; right: respective schemes with artificial viscosity.

### 4.6. Noh problem (2D case)

Finally we consider the Noh test problem [21] in two-dimensional formulation, as it was described in [22]. Here the initial conditions of a polytropic gas with $\gamma = 5/3$ are the following: $\rho = 1$, $p = 0$ (replaced by $10^{-6}$ in the computations), the velocities are directed toward the origin in the *xy* plane with magnitude 1. The solution is a circularly symmetric shock reflecting from the origin: the shock speed is 1/3; behind the shock $\rho = 16$, $p = 16/3$ and $\mathbf{u} = 0$; ahead of the shock the density is $\left(1 + t \, / \sqrt{x^2 + y^2}\right)$, while velocity and pressure remain the same as initially.

The computational domain is $[0, 1] \times [0, 1]$ in the plane. At the boundaries $x = 0$ and $y = 0$ we use the symmetric boundary conditions, while at the boundaries $x = 1$ and $y = 1$ the values of the exact solution (depending on time and radius) are set.

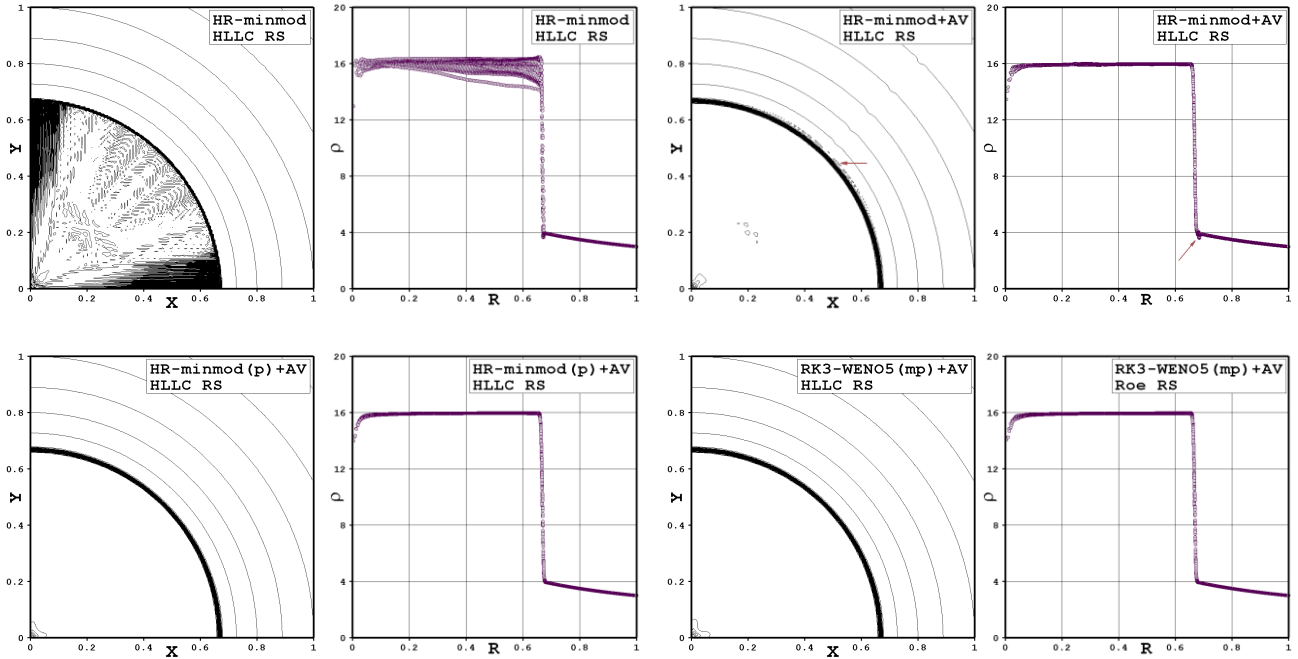

Fig. 17. Results for the Noh problem at *t* = 2 by different schemes. The density contour plots (from 2.5 to 17 with step 0.25) and the scatter plots of density versus radius.

Results for the Noh problem on a grid composed of *J*×*I* = 200×200 cells are shown in Fig. 17. One can see that the HR-minmod scheme with the HLLC solver (top left plots) produces strong shock instability. With adding the artificial viscosity (top right plots) the perturbations behind the shock virtually disappear, although one can notice a small undershoot in density just ahead of the shock (see arrows on plots). This is due to a unique peculiarity of the problem: the flow ahead of the shock is converging and hypersonic. In this specific case the primitive variables are preferred to the characteristic variables when doing the data reconstruction; the left bottom plots in Fig. 17 affirm this statement (the symbol (p) on the plots means using the primitive variables). Based on this finding we propose an alternative version of Suggestion 2: simultaneously with switching to the minmod limiter in the vicinity of the shock one may use the primitive variables there. In other test problems this version works nearly as well as the basic one, but in the current test its advantage is noticeable. The last example in Fig. 17 (bottom right plots) pertains to the RK3-WENO5 scheme with the Roe solver, when employing the artificial viscosity and the alternative version of Suggestion 2.

## 5. Conclusions

The artificial viscosity approach for curing the carbuncle phenomenon was tested on high-order Godunov-type schemes. In so doing, a few selected schemes (HR-minmod, HR-vanLeer, HR-MC, RK2-MC and RK3-WENO5) were used for computing the Quirk-type test problems and other popular tests. The basic theses of the detailed study can be formulated as follows.

1. As applied to high-order schemes, the artificial viscosity approach demonstrates its efficiency. For the principal coefficient of the artificial viscosity model we have chosen the value of $C_{AV} = 0.5$ that ensures the suppression of shock instability if two recommendations to the schemes (presented below) are followed. This value can be considered as an all-purpose one, since it is suitable for computations of any dimensionality and for gases with a wide range of the ratio of specific heats.

2. The first recommendation (Suggestion 1) relates to the procedure of data reconstruction; it is essential to use the characteristic variables here. The use of the primitive variables results frequently in conspicuous oscillations behind the shock. The exception is the reconstruction based on the minmod limiter for which any (primitive or characteristic) variables can be used with similar outputs.

3. The second recommendation (Suggestion 2) refers to the data reconstruction in the cells that get inside the shock layer. One should rest upon the concept of piecewise linear distribution of the data within such cells and use the minmod limiter. In doing so one may apply either the characteristic variables (basic version) or the primitive variables (alternative version). Both variants demonstrated similar results in the majority of test problems, however, as applied to the Noh problem, the alternative version performs better than the basic one.

4. When coupled with the stated recommendations, the artificial viscosity approach not only cures the carbuncle phenomenon, but also reduces substantially the post-shock oscillations. One of the consequences of applying the proposed approach, which is of practical interest, consists in resolving the heating problem that appears when simulating the viscous hypersonic flow around a blunt body.

Our study has also revealed that applying the incomplete Riemann solvers (thought to be free of the carbuncle instability) to the high-order Godunov-type schemes not always ensures the suppression of the carbuncle flaw. Such a situation occurs, for instance, when using the HLL solver in simulations of the Quirk problem on parallelogram grids with $h_y / h_x \leq 1/4$.

Finally, we emphasize that the artificial viscosity approach can be used in a CFD code without its switching or tuning to the specific problem: in a shock-free problem it will have no impact on the solution; in a problem with shock(s) it will provide the necessary dissipation in the shock layer, preventing the carbuncle instability and improving the post-shock solution.

**Acknowledgments**

This work has been supported by the Russian Foundation of Basic Research (Grant 15-01-06224). The author would like to thank Tatiana Zezyulina for the qualified assistance in English and the reviewers for their helpful comments.

## Appendix A. Vortex transport by uniform flow

This appendix provides a convergence study of the selected schemes when simulating the two-dimensional vortex evolution problem for the Euler equations [23].

For this problem we use the computational domain $[0, 10] \times [0, 10]$ with periodic boundary conditions. The initial flow field of a polytropic gas with $\gamma = 1.4$ in the domain is a superposition of two flow patterns: a mean diagonal flow with $(u_x, u_y, \rho, p) = (1, 1, 1, 1)$ and an isentropic vortex (perturbation of the mean flow) defined as

$$\left(\delta u_x, \delta u_y\right) = \frac{\varepsilon}{2\pi} e^{0.5(1-r^2)} \left(-\overline{y}, \overline{x}\right), \qquad \delta T \equiv \delta\left(\frac{p}{\rho}\right) = \frac{(\gamma-1)\varepsilon^2}{8\gamma\pi^2} e^{1-r^2}, \qquad \delta S \equiv \delta\left(\frac{p}{\rho^\gamma}\right) = 0,$$

where $\left(\overline{x}, \overline{y}\right) = \left(x-5, y-5\right)$, $r^2 = \overline{x}^2 + \overline{y}^2$ and the vortex strength $\varepsilon = 5$.

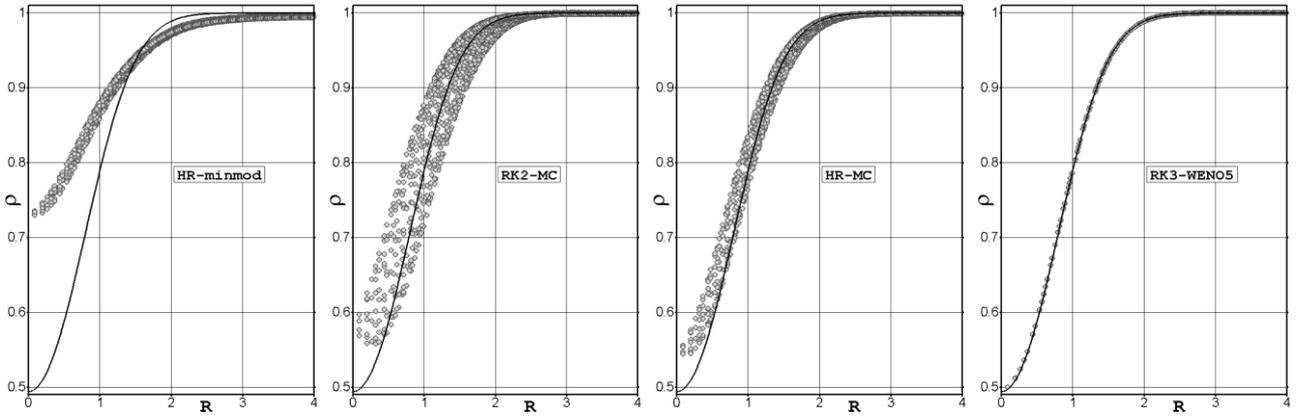

Fig. A.18. Results for the vortex transport problem at $t = 50$ by different schemes with $h = 1/8$: the scatter plots of $\rho$ versus $r$ (solid line depicts the exact solution).

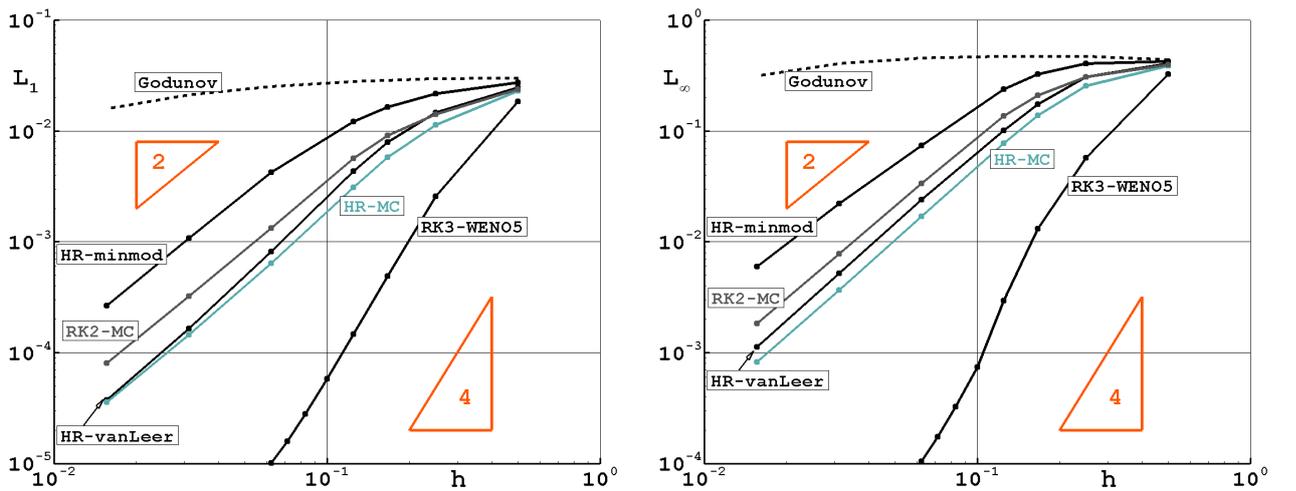

Fig. A.19. Vortex transport problem at t = 50: convergence of density errors with $h$.

The computations were carried out until $t = 50$ (5 periods in time) with $C_{cfl} = 0.6$ for RK3-WENO5 and $C_{cfl} = 0.8$ for the other schemes, using the characteristic variables in all reconstructions. Here we present the data for the selected schemes with the HLLC solver, but it should be noted that sample computations with the exact and Roe solvers gave very close results.

Fig. A.18 shows the scatter plots of density versus distance from the vortex center ($r$) in the case of using the mesh size $h = 1/8$. Fig. A.19 represents the grid convergence study for density errors. As can be seen from this study, the second-order schemes (HR-minmod, HR-vanLeer, HR-MC and RK2-MC) achieve their formal order of accuracy; among them the HR-MC scheme demonstrates the best results. The RK3-WENO5 scheme is much more efficient than the second-order schemes, revealing the convergence order between four and five.

In addition to the high-order schemes, Fig. A.19 displays data obtained by the Godunov scheme. This data indicate that the first-order schemes fail completely to solve the problems of this type.

## References

[1] K.M. Peery, S.T. Imlay, Blunt body flow simulations, AIAA Paper 88-2924, 1988.

[2] J. J. Quirk, A contribution to the great Riemann solver debate, ICASE Report 92-64, 1992; Int. J. Numer. Meth. Fluids 18 (1994) 555–574.

[3] A.V. Rodionov, Artificial viscosity in Godunov-type schemes to cure the carbuncle phenomenon, J. Comput. Phys. 345 (2017) 308–329.

[4] A.V. Rodionov, Artificial viscosity to cure the carbuncle phenomenon: the three-dimensional case, J. Comput. Phys. 361 (2018) 50–55.

[5] B. van Leer, Upwind and high-resolution methods for compressible flow: from donor cell to residual-distribution schemes, Commun. Comput. Phys. 1(2) (2006) 192-206.

[6] R.J. LeVeque, Finite volume methods for hyperbolic problems. Cambridge University Press, 2002.

[7] G.-S. Jiang, C.-W. Shu, Efficient implementation of weighted ENO schemes, J. Comput. Phys. 126 (1996) 202–228.

[8] A. Suresh, H. T. Huynh, Accurate monotonicity-preserving schemes with Runge–Kutta time stepping, J. Comput. Phys. 136 (1997) 83–99.

[9] C.-W. Shu, S. Osher, Efficient implementation of essentially non-oscillatory shock-capturing schemes, J. Comput. Phys. 77 (1988) 439–471.

[10] A. V. Rodionov, Monotonic scheme of the second order of approximation for the continuous calculation of non-equilibrium flows, USSR Comput. Math. Math. Phys. 27 (2) (1987) 175–180.

[11] A. V. Rodionov, Methods of increasing the accuracy in Godunov's scheme, USSR Comput. Math. Math. Phys. 27 (6) (1987) 164–169.

[12] G.D. van Albada, B. van Leer, W.W. Roberts, A comparative study of computational methods in cosmic gas dynamics, Astron. Astrophysics, 108 (1982) 76–84.

[13] B. van Leer, On the relation between the upwind-differencing schemes of Godunov, Engquist–Osher and Roe, SIAM J. Sci. Stat. Comput. 5 (1) (1984) 1–20.

[14] S.K. Godunov, Finite difference method for numerical computation of discontinuous solutions of the equations of fluid dynamics, Mat. Sb. 47 (1959) 271–306.

[15] P.L. Roe, Approximate Riemann Solvers, J. Comput. Phys. 43 (1981) 357–372.

[16] E.F. Toro, M. Spruce, W. Speares, Restoration of the contact surface in the HLL-Riemann solver, Shock Waves 4 (1994) 25–34.

[17] P.R. Woodward, P. Colella, The numerical simulation of two-dimensional fluid flow with strong shocks, J. Comput. Phys. 54 (1984) 115–173.

[18] A. Harten, P.D. Lax, B. van Leer, On upstream differencing and Godunov-type schemes for hyperbolic conservation laws, SIAM Rev. 25 (1) (1983) 35–61.

[19] K. Kitamura, P. Roe, F. Ismail, An evaluation of Euler fluxes for hypersonic flow computations, AIAA Journal 47 (1) (2009) 44–53.

[20] K. Kitamura, E. Shima, Towards shock-stable and accurate hypersonic heating computations: A new pressure flux for AUSM-family schemes, J. Comput. Phys. 245 (2013) 62–83.

[21] W. F. Noh, Errors for calculations of strong shocks using an artificial viscosity and an artificial heat flux, J. Comput. Phys. 72 (1987) 78–120.

[22] R. Liska, B. Wendroff, Comparison of several difference schemes on 1D and 2D test problems for the Euler equations, SIAM J. Sci. Comput. 25 (3) (2003) 995–1017.

[23] C. Hu, C.-W. Shu, Weighted essentially non-oscillatory schemes on triangular meshes, J. Comput. Phys. 150 (1999) 97–127.